\documentclass[a4paper,11pt]{article}
\usepackage[multidot]{grffile}
\usepackage{latexsym,amssymb,enumerate,amsmath,epsfig,amsthm}
\usepackage[margin=1in]{geometry}
\usepackage{setspace,color}
\usepackage{graphics}
\usepackage{graphicx}
\usepackage[ruled]{algorithm2e}
\usepackage{empheq}
\usepackage{bm,amsmath}
\usepackage{subcaption}
\usepackage{hyperref} 
\hypersetup{
    colorlinks=true,
    linkcolor=blue,
    filecolor=magenta,      
    urlcolor=cyan,
    citecolor=red,
    pdfpagemode=FullScreen,
    }

\usepackage{multirow}
\newcommand\fixup{\kern-\fontcharic\scriptfont2`\"}

\newcommand{\p}{\mathbf{p}}
\newcommand{\q}{\mathbf{q}}

\title{A Semi-Lagrangian Spherical Essentially Non-Oscillatory (SENO) Scheme for Advection Equations of $\mathbb{S}^2$-valued Functions}
\author{
Shingyu Leung\thanks{Department of Mathematics, the Hong Kong University of Science and Technology, Clear Water Bay, Hong Kong. Email: {\bf masyleung@ust.hk}}
}

\date{}

\begin{document}
\thispagestyle{plain}
\maketitle

\begin{abstract}
We develop a numerical scheme for solving the advection equation of $\mathbb{S}^2$-valued functions of real variables, which models the time-evolution of a $\mathbb{S}^2$-valued mapping on the real line by a known velocity field. The idea is to extend the semi-Lagrangian method for the linear scalar advection equation. We first construct the backward flow map between two adjacent time levels and then interpolate the discrete ordered data of $\mathbb{S}^2$. To handle $\mathbb{S}^2$-functions which have kinks or sharp discontinuity in their components, we incorporate the \textit{Spherical Essentially Non-Oscillatory} (SENO) interpolation method, which effectively reduces the spurious oscillations in high-order reconstructions. We will show multiple examples to demonstrate the accuracy and effectiveness of the proposed algorithm for the partial differential equation of $\mathbb{S}^2$-functions.
\end{abstract}

\section{Introduction}
\label{Sec:Introduction}

This paper constructs high-order numerical schemes for solving the advection equation of $\mathbb{S}^2$-functions of real variables, given by
\begin{equation}
\p_t + c(s,t) \p_s = 0 
\label{Eqn:advection}
\end{equation}
with an initial condition $\p(s,0)=\p_0(s)$ for $\mathbb{S}^2$-valued functions $\p=\p(s,t): [0,1]\times[0,\infty) \rightarrow \mathbb{S}^2$ and $c(s,t):[0,1]\times[0,\infty) \rightarrow \mathbb{R}$ representing the speed of propagation. The parameter $s$ represents the parameterization of data points on a unit sphere. For simplicity, we assume a periodic boundary condition such that $\p(0,t)=\p(1,t)$. The variable $t$ is a time-like variable representing the evolution of these data points on the unit sphere $\mathbb{S}^2$. Since we can represent $\p$ in a component form, i.e. $\p(s,t)=(x(s,t),y(s,t),z(s,t))$ with $x,y$ and $z \in \mathbb{R}$ and $(x,y,z)\in\mathbb{S}^2$, we can obtain the analytical solution to the advection equation (\ref{Eqn:advection}) using the method of characteristics. Therefore, we have $\p(s(t),t)=\p_0(s_i)$ along a characteristic $s'(t)=c(s(t),t)$ with the initial take-off location $s=s_i$ at $t=0$. When we interpret the initial condition $\p_0(s)$ as a parameterized closed curve on $\mathbb{S}^2$, the solution to the advection equation $\p(s,t)$ gives a re-parameterization of the same curve at a later time $t$. There are applications related to the interpolation problems in $\mathbb{S}^2$  and involving $\mathbb{S}^2$-valued functions. For example, we can find applications in quantum field theory from quantum mechanics \cite{adl86}, modeling protein structures \cite{pro14}, molecular dynamics simulation \cite{rap85}, {particle dynamics} in fluid mechanics \cite{ghkr06}, fluid flow visualizations \cite{hanma95}, computations of flexible filaments and fibers in complex fluids \cite{tscfro20,stwk21}, and some applications to dynamics of rigid-bodies \cite{weiterfed06,udwsch10}.

Since components in the three-vector are decoupled, one trivial approach to this problem is to directly apply the finite difference upwind scheme and update the solution using
$$
\frac{\p_i^{n+1}-\p_i^n}{\Delta t} + c_i^+ \cdot \frac{\p_i^n-\p_{i-1}^n}{\Delta s} + c_i^- \cdot \frac{\p_{i+1}^n-\p_i^n}{\Delta s} =0
$$
where $c_i^{\pm}=(c(s_i,t^n) \pm |c(s_i,t^n)|)/2$ which involves three independent hyperbolic partial differential equation (PDE). We can also apply high-order methods like TVD-RK and ENO/WENO \cite{shuosh88,shuosh89,liuoshcha94,shu97} to obtain more accurate solutions. However, the main issue of this approach is that there is no guarantee that the solution stays on $\mathbb{S}^2$ in the evolution, which might cause nonphysical applications in practice. 

Therefore, in this work, we will develop a simple numerical scheme to solve the advection equation (\ref{Eqn:advection}) based on the SENO interpolation proposed in \cite{fonleu23}. Extending the method of characteristics for the scalar advection equation, we develop a semi-Lagrangian method for the advection equation of $\mathbb{S}^2$-functions. The first step of the algorithm relies on the accurate construction of the backward flow map, which identifies the takeoff location between two time levels under the advection field. Since the function is constant along this characteristic, we can determine the solution at the later time level by interpolating the function defined at the discretized sampling location. Now, because we are looking at $\mathbb{S}^2$-valued functions, these discrete data give a set of ordered points on a unit sphere.

Several methods exist to interpolate these data points on the unit sphere. For example, we have the spherical linear interpolation (SLERP), and the spherical quadrangle interpolation (SQUAD) based on the quaternion representation \cite{shoemake_85}. These two are the most popular and commonly used interpolation methods on the unit sphere. The SLERP interpolation provides the piecewise linear interpolation in geodesic on the unit sphere, while the SQUAD gives a smooth and slightly higher-order reconstruction of the data points. In \cite{fonleu23}, we have developed an interpolation scheme named the \textit{Spherical Interpolation of orDER} $n$ (SIDER-$n$) that produces a $C^{n}$ interpolant given $n \geq 2$. The idea follows the construction of the B\'{e}zier curves based on the composition of multiple SLERPs. We have also followed the ENO philosophy and have proposed a new \textit{Spherical Essentially Non-Oscillatory} (SENO) interpolation method. This approach can provide a smooth, high-order interpolant even when the underlying curve has kinks. Therefore, we propose incorporating this SENO idea into the backward flow map and obtaining a high-order semi-Lagrangian method for the advection equation.

The rest of the paper is organized as follows. We will summarize the proposed spherical interpolation of order-$n$ (SIDER-n) scheme and the spherical essentially non-oscillatory (SENO) scheme as developed in \cite{fonleu23} in Section \ref{Sec:SENO}. Section \ref{Sec:Proposal} gives our proposed algorithm for the advection equation of $\mathbb{S}^2$-functions based on the semi-Lagrangian scheme. Some examples will be given in Section \ref{Sec:Examples} to demonstrate the accuracy of our proposed numerical approaches.

\section{The Spherical Essentially Non-Oscillatory (SENO) Scheme}
\label{Sec:SENO}

Quaternions are numbers consisting of four dimensions, one real part and a three-dimensional analogy to {the imaginary part of} complex numbers \cite{ham63}. A quaternion can be written in many forms: \[\underset{\text{real}}{\boxed{a}} + \underset{\text{imaginary}}{\boxed{b \mathbf{i} + c\mathbf{j} + d\mathbf{k}}} = (a, b, c, d) = (\underset{\text{scalar}}{\boxed{a}}, \underset{\text{vector}}{\boxed{\mathbf{u}}}),\] where \(a, b, c, d \in \mathbb{R}\), $\mathbf{u} = (b, c, d) \in \mathbb{R}^3$. The notations $\mathbf{i}$, $\mathbf{j}$ and $\mathbf{k}$ are extensions of {the imaginary part of} complex numbers with the properties that $\mathbf{i}^2=\mathbf{j}^2=\mathbf{k}^2{=\mathbf{ijk}}=-1$, $\mathbf{ij}=\mathbf{k}$ with the {bicyclic permutation with respect to $\mathbf{i}$ that $1\rightarrow \mathbf{i}\rightarrow -1\rightarrow -\mathbf{i}$ and $\mathbf{j}\rightarrow -\mathbf{k}\rightarrow -\mathbf{j}\rightarrow \mathbf{k}$}. Some important properties about quaternions include \begin{itemize}
\item Hamilton product. $(a_1, \mathbf{u_1})(a_2, \mathbf{u_2}) = (a_1a_2 - \mathbf{u_1} \cdot \mathbf{u_2}, a_1\mathbf{u_2} + a_2\mathbf{u_1} + \mathbf{u_1} \times \mathbf{u_2})$,
where the notation {$\cdot$ and $\times$ denotes the typical dot and cross product.}
\item Inverse map. $\mathbf{q}^{-1} = (a, -\mathbf{u})/(a^2 + b^2 + c^2 + d^2)$. If $\mathbf{q} = (a, \mathbf{u})$. In particular, if $\mathbf{q}$ is a unit quaternion, $\mathbf{q}^{-1} = (a, -\mathbf{u})$.
\item Exponential map. $\exp(a, \mathbf{u}) = \exp(a)(\cos \lVert\mathbf{u}\rVert, ((\sin \lVert\mathbf{u}\rVert)/\lVert\mathbf{u}\rVert) \mathbf{u})$ {where the norm notation $\lVert \cdot \rVert$ denotes the 2-norm in this paper, unless otherwise specified.}
\item Logarithm map. $\ln(a, \mathbf{u}) = \left(\ln \sqrt{a^2 + \lVert \mathbf{u} \rVert^2}, \frac{1}{\lVert \mathbf{u} \rVert} 
\arccos \left( \frac{a}{\sqrt{a^2 + \lVert \mathbf{u} \rVert^2}} \right) 
\mathbf{u} \right)$. 
\item Power map 
\begin{align*}
&\hspace{1.33em} (a, \mathbf{u})^{f(t)} = \exp(f(t)\ln(a, \mathbf{u})) = \exp(f(t)\ln \sqrt{a^2 + \lVert \mathbf{u} \rVert^2}, (f(t)k/\lVert \mathbf{u} \rVert) \mathbf{u}) \\
&= \left((a^2 + \lVert \mathbf{u} \rVert^2)^{f(t)/2}\cos (f(t)k), (a^2 + \lVert \mathbf{u} \rVert^2)^{f(t)/2} \left[\sin (f(t)k)/\lVert\mathbf{u}\rVert\right] \mathbf{u}\right) \, ,
\end{align*}
where $k = \arccos \left(a/\sqrt{a^2 + \lVert \mathbf{u} \rVert^2}\right)$. When a quaternion has its {2-}norm equal to one, we call them a unit quaternion. If $(a, \mathbf{u})$ is a unit quaternion, then 
$$
(a, \mathbf{u})^{f(t)} = \left(\cos (f(t)k), \left[ \frac{\sin (f(t)k)}{\lVert \mathbf{u} \rVert} \right] \mathbf{u} \right) \, .
$$
\end{itemize}

Because we can use these unit quaternions to define rotation, we also call these quaternions rotation quaternions. With proper definitions, they can rotate a position vector defined in either $\mathbb{S}^2$ or $\mathbb{R}^3$ while preserving the length of the vector. To see this, we express a unit quaternion as $(a, b, c, d) = (a, \mathbf{u}) = (\cos (\theta/2), {\sin (\theta/2)}\mathbf{v})$ where $\mathbf{v}$ is a unit vector representing the 3D rotation axis, and $\theta$ is the anticlockwise{/counterclockwise} rotation angle around $\mathbf{v}$ carried by the rotation quaternion. If we want to rotate $\mathbf{p_a}\in \mathbb{S}^2$ with a rotation quaternion $\mathbf{r_{ab}} = (\cos (\theta_{ab}/2), \sin (\theta_{ab}/2)\mathbf{a_{ab}})$, we can first convert $\mathbf{\mathbf{p_a}}$ to a unit quaternion given by $\mathbf{\mathbf{q_a}} = (0, \mathbf{\mathbf{p_a}})$. Then we apply the rotation operator given by
$\mbox{ROTATE}(\mathbf{q_a}, \mathbf{r_{ab}}) = (\mathbf{r_{ab}})(\mathbf{q_a})(\mathbf{r_{ab}})^{-1}$. The final position after the rotation is given by the {imaginary part of the} unit quaternion $\mathbf{q_b} = (0, \mathbf{p_b})$. 

Introducing a parameterization $t$ such that $t=0$ and $t=1$ corresponding to the initial position $\mathbf{q_a}$ and $\mathbf{q_b}$, respectively, we can interpolate these two data points by the rotation operator $\mbox{ROTATE}(\mathbf{q_a}, \mathbf{r_{ab}}, t) = (\mathbf{r_{ab}})^t (\mathbf{q_a}) (\mathbf{r_{ab}})^{-t}$ for $t \in [0, 1]$. This expression leads to the so-called SLERP (\textit{Spherical Linear intERPolation}) formula \cite{shoemake_85,sola}: 
$\mbox{SLERP}(\mathbf{\mathbf{q_a}}, \mathbf{q_b}, t) = (\mathbf{q_a}) ((\mathbf{q_a})^{-1} \mathbf{q_b})^t$ for $t \in [0, 1]$.
In particular, if the quantity $\mathbf{a_{ab}}$ in the rotation quaternion is perpendicular (and this is assumed to be true in the remaining of this article) to both $\mathbf{\mathbf{p_a}}$ and $\mathbf{p_b}$, we have $\mathbf{\mathbf{p_a}} \cdot \mathbf{p_b} = \cos \theta_{ab}$ and 
\begin{equation*}
    \begin{split}
        \mathbf{p_b} &= (\cos \theta_{ab}) \mathbf{\mathbf{p_a}} + (\sin \theta_{ab}) (\mathbf{a_{ab}} \times \mathbf{\mathbf{p_a}}); \quad \quad \text{(Rodrigues' rotation formula \cite{rodrigues_40})}\\
        \mathbf{a_{ab}} &= \mathbf{\mathbf{p_a}} \times \left[\dfrac{\mathbf{p_b} - (\cos \theta_{ab}) \mathbf{\mathbf{p_a}}}{(\sin \theta_{ab})}\right] = \dfrac{\mathbf{\mathbf{p_a}} \times \mathbf{p_b}}{\sin \theta_{ab}}.
    \end{split}
\end{equation*}

In the recent paper \cite{fonleu23}, we have introduced a new class of interpolation schemes on the unit sphere, denoted by SIDER. With reference to the construction of quadratic B\'{e}zier curves, we have developed the following spherical quadratic curve (denoted by SIDER2),
\begin{eqnarray*}
& &\mbox{SIDER2}(\boxed{\underset{\text{start}}{\mathbf{q_1}}}, \boxed{\underset{\text{second data}}{\mathbf{q_2}}}, \boxed{\underset{\text{end}}{\mathbf{q_3}}}, t) \\
&=&\mbox{SLERP}(\mbox{SLERP}(\mathbf{q_1},\mathbf{\mathbf{d_{2a}}}, t), \mbox{SLERP}(\mathbf{\mathbf{d_{2b}}}, \mathbf{q_3},t), t)\\
&=&\mathbf{q_1} \left(\mathbf{q_1}^{-1} {\mathbf{d_{2a}}}\right)^t \left[ \left[\mathbf{q_1}(\mathbf{q_1}^{-1} \mathbf{d_{2a}})^{t} \right]^{-1} \left[\mathbf{d_{2b}}(\mathbf{d_{2b}}^{-1} \mathbf{q_3})\right]^t\right]^{f_2(t)} {= 
\left(0,\p_{\mbox{SIDER2}}(t) \right)}
\end{eqnarray*}
where $t \in [t_1, t_3]$, and $f_2(t) = (t - t_1)/(t_3 - t_1)$. Unless specified otherwise, we might use $t_1 = 0$ and $t_3 = 1$. The points $\mathbf{q_i} = (0, \mathbf{\mathbf{p_i}})$, $\mathbf{d_{2a}} = (0, \mathbf{\mathbf{c_{2a}}})$ and $\mathbf{d_{2b}} = (0, \mathbf{\mathbf{c_{2b}}})$ are the quaternion representation of the position vectors $\mathbf{\mathbf{p_i}}$, $\mathbf{\mathbf{c_{2a}}}$ and $\mathbf{\mathbf{c_{2b}}}$, respectively. We construct $\mathbf{c_{2b}}$ (and $\mathbf{c_{2a}}$) using the geodesic extrapolating based on the first data points $\mathbf{p_1}$ (and $\mathbf{p_3}$) and the intermediate one $\mathbf{p_2}$ so that the final interpolant reaches $\mathbf{p_2}$ when $t = t_2 = 0.5(t_1 + t_3)$. To enforce this condition, we refer to the spatial relationships among the data points and the (only) control point in a quadratic force interpolating B\'{e}zier curve. Mathematically, we assign $\mathbf{d_{2a}} = (0, \mathbf{\mathbf{c_{2a}}}) = \text{SLERP}(\mathbf{q_3},\mathbf{q_2},2)$ and $\mathbf{d_{2b}} = (0, \mathbf{\mathbf{c_{2b}}}) = \text{SLERP}(\mathbf{q_1}, \mathbf{q_2},2)$.

Similarly, we have the following higher order extension
\begin{eqnarray*}
&& \mbox{SIDER3}(\boxed{\underset{\text{start}}{\mathbf{q_1}}}, \boxed{\underset{\text{second point}}{\mathbf{q_2}}}, \boxed{\underset{\text{third point}}{\mathbf{q_3}}}, \boxed{\underset{\text{end}}{\mathbf{q_4}}}, t) \\
&=& \mbox{SLERP}(\mbox{SIDER2}(\mathbf{q_1}, \mathbf{q_2},\mathbf{q_3},g_3(t)), \mbox{SIDER2}(\mathbf{q_2},\mathbf{q_3},\mathbf{q_4},h_3(t)), f_3(t)) \\
&=& \left(0,\p_{\mbox{SIDER3}}(t) \right)
\end{eqnarray*}
where $t \in [t_1, t_4]$. Therefore, a SIDER3 reconstruction is a linear combination of two scaled SIDER2, that we interpolate within $\{\mathbf{p_1}, \mathbf{p_2}, \mathbf{p_3}\}$ and $\{\mathbf{p_2}, \mathbf{p_3}, \mathbf{p_4}\}$ simultaneously. For simplicity, we set the starting time $t_1 = 0$ and the ending time $t_4 = 1$, so that $t_2 = \frac{1}{3}$ and $t_3 = \frac{2}{3}$, and $g_3(t) = 3t/2$, $h_3(t) = (3t - 1)/2$ and $f_3(t) = t$. 

One usually observes oscillations in the interpolant when reconstructing a high-order curve with sharp changes and turns, and this behavior is undesirable in many applications. In \cite{fonleu23}, we follow the philosophy of Essentially Non-Oscillatory (ENO) and propose an ENO interpolation on the unit sphere. We have named the interpolation approach the \textit{Spherical Essentially Non-Oscillatory} (SENO in short). Given a set of $2n$ data points, denoted by $\mathbf{p_{i-n+1}}, \cdots, \mathbf{p_i}$, $\mathbf{p_{i+1}}, \cdots \mathbf{p_{i+n}}$, we are interested in constructing a high-order curve between $\mathbf{p_i}$ and $\mathbf{p_{i+1}}$. To do this, we first reconstruct a $C^{n}$ curve from any $n+1$ consecutive data points using SIDER-$n$. For example, for $n=2$, i.e., we are given four data points, we first construct two $C^2$ SIDER2 curves from any three consecutive data on the unit sphere. When $n=3$, we have in total six data points. From these, we obtain three $C^3$ curves obtained by SIDER3. To avoid an oscillatory interpolant, we consider these $n$ interpolants from {SIDER}-$(n-1)$ and determine the corresponding \textit{variation} of these curves between the data points $\mathbf{p_{i}}$ and $\mathbf{p_{i+1}}$. The one with the least \textit{variation} is chosen to represent the SENO interpolant between the points $\mathbf{p_i}$ and $\mathbf{p_{i+1}}$.

\section{Our Proposed Approach}
\label{Sec:Proposal}

\subsection{The Semi-Lagrangian Scheme}

Let $\Psi_a^b(\mathbf{y}):\mathbb{R}^d\rightarrow \mathbb{R}^d$ be the flow map associated with the ordinary differential equation $\mathbf{x}'(t)=c(\mathbf{x},t)$ with the initial condition $\mathbf{x}(a)=\mathbf{y}$ so that $\Psi_a^b(\mathbf{y})=\mathbf{x}(b)$. Motivated by the method of the characteristics, the semi-Lagrangian scheme for the advection equation of a scalar function 
$f_t+c(s,t) f_s=0$ is given by $f(s_i,t^{n+1}) \approx f_i^{n+1} = f\left( \Psi_{t^{n+1}}^{t^n}(s_i), t^n \right)$.
This scheme involves two main steps, including the construction of the backward flow map $\Psi_{t^{n+1}}^{t^n}$ and the interpolation step involving the sampled solution at the time level $t=t^n$. For $\mathbb{S}^2$-valued functions such as (\ref{Eqn:advection}), we notice that the unknown function is vectored-valued with the domain staying in the Cartesian space. Therefore, we can simply keep the flow map construction stage while replacing the typical one-dimensional interpolation of $f$ with an interpolation scheme for $\p\in \mathbb{S}^2$. In particular, we have $\p_i^{n+1} = \p \left( \Psi_{t^{n+1}}^{t^n}(s_i), t^n \right)$. {This approach is sometimes referred to as the \text{strong form} of the semi-Lagrangian method.} The first step of the algorithm is to construct the backward flow map from $t=t^{n+1}$ to $t=t^n$ by solving the dynamical system $s'(t)=c(s,t)$ with the terminal condition $s(t^{n+1})=s_i$ for $s(t^n)$. This step can be done using typical high-order numerical methods such as Runge-Kutta schemes. The main issue is the second step involving the interpolation problem. In particular, we are given a set of sampled values $\{ (s_i, \p_i^n) \}$ on a set of uniform mesh $\{s_i\}$ and are required to interpolate the value of $\p$ at a generally non-mesh point $s= \Psi_{t^{n+1}}^{t^n}(s_i)$.

\subsection{Component-wise Interpolations}

Since the equation is decoupled, one might perform the interpolation problem in a component-by-component fashion. Since the interpolation problem is the foundation of many advanced numerical algorithms, many available implementations exist. For example, \textsf{MATLAB} has a simple piecewise linear interpolation (the default implementation in \textsf{interp1}). When one needs high-order accurate solutions, we can consider using \textsf{pchip} representing the monotone high-order shape-preserving piecewise cubic interpolation. Although the scheme sounds straightforward, the algorithm does not consider that the resulting solution has to be on $\mathbb{S}^2$. There is no guarantee (and there is no reason we would expect) that these interpolation methods produce an interpolant that stays on the unit sphere. In particular, the piecewise linear interpolation simply joins two adjacent data points on $\mathbb{S}^2$ using a straight line in $\mathbb{R}^3$ but does not lead to the geodesic curve on $\mathbb{S}^2$. 

\begin{figure}[!h]
\centering
(a)\includegraphics[width=0.32\textwidth]{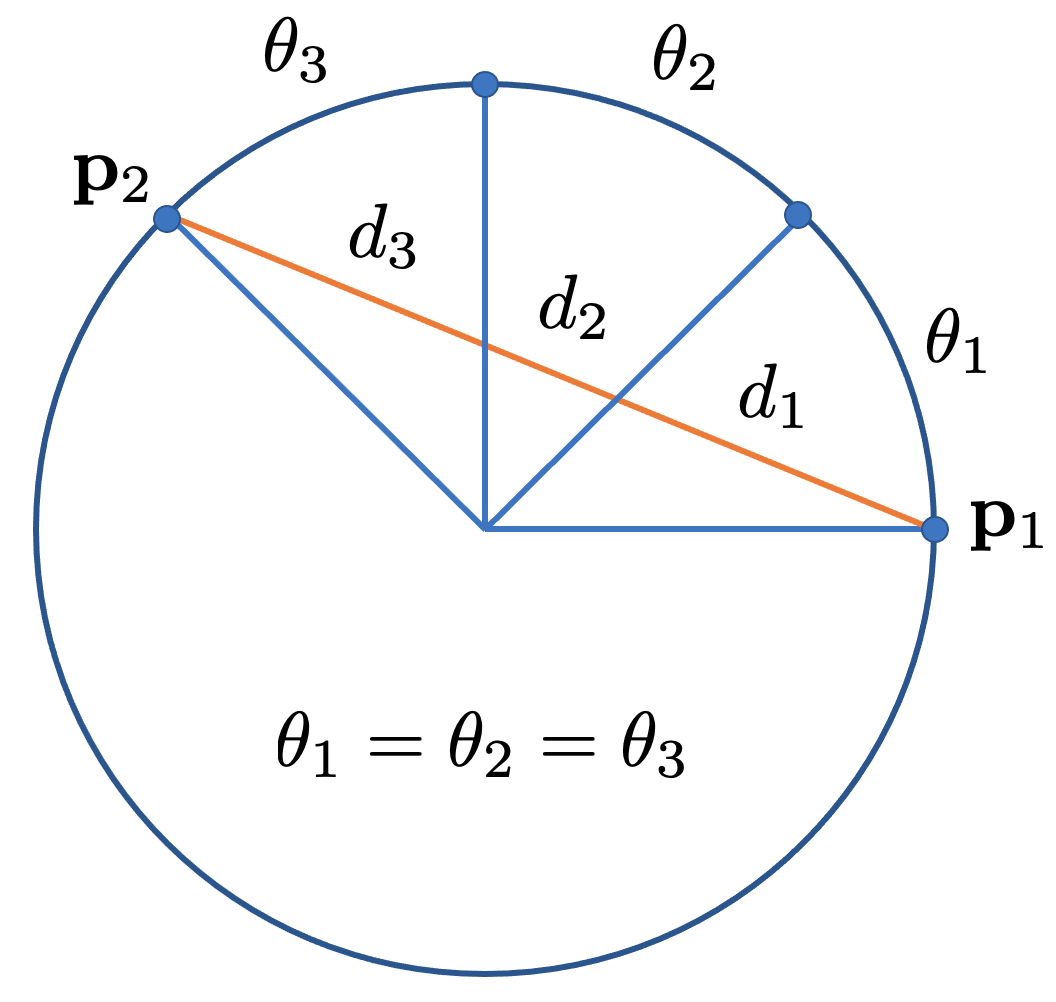} 
(b)\includegraphics[width=0.32\textwidth]{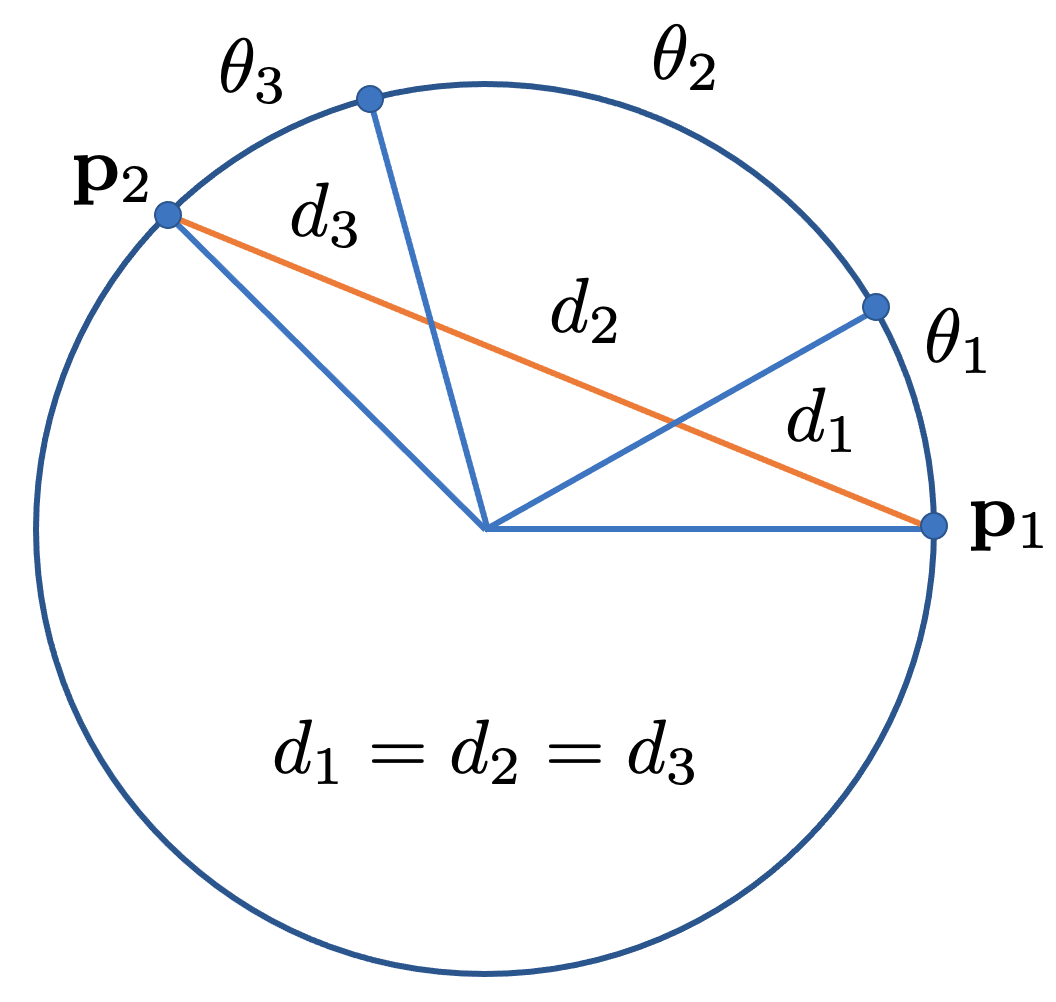} 
\caption{(a) Uniform partition on the unit sphere and (b) uniform partition on the straight line joining two points on a sphere.} 
\label{Fig:UniformPartition}
\end{figure}

One possible resolution is to incorporate a projection step after each interpolation, i.e.
$\p(s_i,t^{n+1}) =\p \left( \Psi_{t^{n+1}}^{t^n}(s_i), t^n \right) / 
\|\p \left( \Psi_{t^{n+1}}^{t^n}(s_i), t^n \right)\|$
so that the solution for each $s_i$ has a unit length. However, such a projection step does not preserve a uniform sampling. For example, when we uniformly partition the straight line joining two non-opposite points on the unit sphere in $\mathbb{R}^3$, the data points projected onto the unit sphere will generally not be uniformly separated. We have shown a simple demonstration in Figure \ref{Fig:UniformPartition} where we plot the two data points $\p_1$ and $\p_2$ on the great circle. The geodesic distance between these two points is determined by the angle $\theta$, while the Euclidean distance is given by $2 \sin \frac{\theta}{2}$. Now, if we construct $n$-partitions so that each segment has the length $\frac{2}{n} \sin \frac{\theta}{2}$ and then project the first point away from the sphere onto the sphere, the geodesic distance from this projection to the endpoint of the interval (denoted by $\beta$) satisfies
${\sin\beta}/{\sin (\frac{\theta}{2}-\beta)} = \frac{2}{n} \sin \frac{\theta}{2}$. Since the solution is different from $\beta=\theta/n$ for general $n$ and $\theta$, a uniform partition along the straight line joining two adjacent data points will provide projected sampling locations closer to the given data points $\p_1$ and $\p_2$. This observation also implies that the semi-Lagrangian method based on componentwise linear interpolation prefers to generate positions near the given sampling points.

\subsection{The Semi-Lagrangian Scheme with the SENO Interpolation}

To respect the geometry of the interpolation problem, we propose to replace all componentwise interpolation methods with the SLERP, our proposed spherical interpolation schemes SIDER, or the non-oscillatory schemes SENO. We discuss in detail the approach based on SENO3, and the implementation based on SLERP or SENO2 is relatively straightforward.

Assuming that we are given the discretized solution $\p_i^n$ on a uniform mesh $s_i\in[0,1]$ at the \textit{time} level $t=t^n$, we first construct the backward flow map $\Psi_{t^{n+1}}^{t^n}(s_i)$ by solving the corresponding ODE backward in time using any well-developed numerical integrator. Then for each grid point $s_i$, we determine the index $j$ such that $\Psi_{t^{n+1}}^{t^n}(s_i)\in [s_j,s_{j+1}]$. According to the SENO3 construction, we first obtain the following three SIDER3 interpolants,
\begin{eqnarray*}
&& \mbox{SIDER3}\left( \q_{j-2},\q_{j-1},\q_j,\q_{j+1}, r \right) \, \mbox{ for $\frac{2}{3}\le r \le 1$}\\
&& \mbox{SIDER3}\left( \q_{j-1},\q_{j},\q_{j+1},\q_{j+2}, r \right) \, \mbox{ for $\frac{1}{3}\le r \le \frac{2}{3}$}\\
&& \mbox{SIDER3}\left( \q_{j},\q_{j+1},\q_{j+2},\q_{j+3}, r \right) \, \mbox{ for $0\le r \le \frac{1}{3}$}
\end{eqnarray*}
where $\q_i=(0,\p_i)$ is the quaternion representation of $\p_i$. Then we evaluate one of these interpolants at either
$$
\frac{2}{3}+ \frac{\Psi_{t^{n+1}}^{t^n}(s_i)-s_j}{3\Delta s} \, , \, \frac{1}{3}+ \frac{\Psi_{t^{n+1}}^{t^n}(s_i)-s_j}{3\Delta s} \, \mbox{ or } \, \frac{\Psi_{t^{n+1}}^{t^n}(s_i)-s_j}{3\Delta s}
$$
depending on which set of stencils gives an interpolant that produces the shortest geodesic distance on the corresponding interval.

For the rest of this section, we discuss properties of the method and give some implementation details. 

\paragraph{Stability Condition.} Note that when using the simple first-order Euler method scheme to construct the flow map with the time step $\Delta s$ satisfying the condition $\Delta s< \Delta x/\max|c(s_i,t^n)|$, the numerical method is equivalent to the typical finite difference upwind scheme. However, unlike the finite difference method, the semi-Lagrangian method does not require the above inequality condition to assure the algorithm's stability. Because of the search of the interval $[s_j,s_{j+1}]$ that contains the takeoff location, one can prove the required stability in the numerical solution at the $t=t^{n+1}$ level. Having said that, it does not imply that the timestep can be arbitrarily large. The stability constraint for $\Delta s$ in the semi-Lagrangian comes from the absolute stability condition of the chosen numerical integrator for the ODE, which usually depends on $|c'(s,t)|$. 

\paragraph{Guarantee of a $\mathbb{S}^2$-Solution.} The backward flow map obtained in the first step of the proposed approach only relates the spatial locations $s$ at different times. It relies on the second stage of the algorithm to guarantee the constraint that the numerical solution should be data points on $\mathbb{S}^2$. Since each SENO interpolation consists of multiple calls of SLERP reconstruction and SLERP guarantees to give data points on the unit sphere, our semi-Lagrangian scheme can automatically produce a $\mathbb{S}^2$-function as the numerical solution to the advection equation (\ref{Eqn:advection}). 

\paragraph{Computational Efficiency.} There are multiple ways to further improve the computational efficiency of the semi-Lagrangian method. For example, if the velocity field is autonomous such that the velocity is independent of time, the flow map is independent of time but depends on the difference $t^{n+1}-t^n$. Mathematically, we have $\Psi_{t^1}^{t^0} = \Psi_{t^2}^{t^1} =\cdots = \Psi_{t^{n+1}}^{t^n}$. We can reuse the flow map for all solution updates. Moreover, if we want to reduce the number of interpolation steps to obtain the final solution at the final time $t=t^{2^k}$, we have $\Psi_{t^{2^k}}^{t^0}= \Psi_{t^{2^{k-1}}}^{t^0} \circ \Psi_{t^{2^k}}^{t^{2^{k-1}}} =  \left( \Psi_{t^{2^{k-1}}}^{t^0} \right)^2 = \left[ \left( \Psi_{t^{2^{k-2}}}^{t^0} \right)^2 \right]^2 =\cdots$. We can recursively define these intermediate flow maps and, therefore, obtain the overall flow map $\Psi_{t^{2^k}}^{t^0}$ in only $k$-iterations. When the velocity field is time-dependent but periodic, we can still solve the ODE for the period $T$ and reuse it for a more extended propagation. We refer interested readers to \cite{canyin06,leuqia10,leu13} and thereafter for more discussions on this flow map doubling approach and efficient implementation of the semi-Lagrangian method. 

\paragraph{An Eulerian Interpolation Scheme.} Some recent work \cite{youwonleu17,youleu18,youleu20} have also developed an Eulerian interpolation scheme to construct the high-order flow map on a uniform mesh based on interpolations. This method can naturally replace the ODE problem in the semi-Lagrangian step with a PDE problem. In particular, if we solve the following advection equation
\textit{forward} in time $\phi_t(s,t)+c(s,t) \phi_s(s,t)=0$ with the initial condition $\phi(s,t^n)=s$ using any well-developed finite difference schemes on a uniform mesh such as TVDRK-WENO methods, the \textit{backward} flow map $\Psi_{t^{n+1}}^{t^n} (s)$ is given by $\phi(s,t^{n+1})$, i.e. $\Psi_{t^{n+1}}^{t^n} (s) = \phi(s,t^{n+1})$. We have discussed this Eulerian framework in detail when developing a level-set method \cite{oshset88,set96,oshfed03} for computational dynamical systems in a series of studies in \cite{leu11,youleu14} and thereafter.

\paragraph{Complexity.} Indeed, it is possible to construct the backward flow map $\Psi_{t^n}^{t^0}$ directly and interpolate the initial condition only once. This observation implies that $\p_i^n=\p_0 \left( \Psi_{t^n}^{t^0}(s_i) \right)$. When the initial condition is explicitly given as a function, we only need to evaluate the initial function of $\mathbb{S}^2$ on a specific set of (usually nonuniform) locations. When the initial condition is given on a set of uniformly sampled points $s_i$, we perform one single interpolation on a set of non-mesh points to obtain the solution to the PDE at the final time $t=t^n$. But this approach does not provide all intermediate solutions for visualization. To get the solution at a specific time level, one can always trace the characteristic back to $t=t^0$ and perform one interpolation. This gives the so-called \textit{global} Lagrangian approach (as discussed in \cite{leuqia10}), which effectively reduces the total number of interpolation steps. But when we want to observe the evolution of the solution in time, the computational complexity could be extremely high. In particular, let $N$ be the number of mesh points and $K$ be the complexity for advancing the solution for a single time step. Since the computational complexity of SLERP and SENO interpolations are all $O(N)$, the overall complexity of the semi-Lagrangian method for determining solutions at all $n$-levels is $O(n(K+N))$, while that of the global Lagrangian method is $O(n^2K+nN)$.

\paragraph{High-dimensional Problems.} High-dimensional generalizations are straightforward. For example, when solving the advection equation 
$\p_{t}(\mathbf{s},t) + \mathbf{c}(\mathbf{s},t) \cdot \nabla_{\mathbf{s}} \p(\mathbf{s},t) =0$, we can first obtain the flow map by solving the corresponding high-dimensional ODE in the Cartesian space and interpolate the $\mathbb{S}^2$-data on the mesh using SENO in a dimension-by-dimension fashion. The idea is the same as solving multi-dimensional scalar hyperbolic equations using high-order ENO/WENO schemes. Instead of designing a multi-dimensional SENO to interpolate the function at once, we split the problem into multiple one-dimensional interpolation problems.

\section{Numerical Examples}
\label{Sec:Examples}

This section considers several numerical examples to confirm the numerical accuracy of the proposed numerical schemes. We have two approaches to construct the initial condition $\p_0(s)$ for $s\in[0,1]$. The first approach defines a function $g(s)$ on the cylindrical surface $x^2+y^2=1$ and then computes the projection onto the unit sphere 
\begin{equation}
\p_0(s)=\left( \frac{\cos 2 \pi s}{\sqrt{1+g^2(s)}} \, , \, \frac{\sin 2 \pi s}{\sqrt{1+g^2(s)}} \, , \,  \frac{g(s)}{\sqrt{1+g^2(s)}} \right) \, .
\label{Eqn:ExProjection1}
\end{equation}
The second initialization defines two curves $z=h_{\pm}(y)$ on the planes $x=\pm 1$, respectively, and then projects them onto $\mathbb{S}^2$. This construction provides a simple way to define a discontinuous condition on the unit sphere. In the first two examples, we consider a smooth initial condition $g_1(s)=\sin(20\pi s)$ and an initial condition with kinks given by $g_2(s)=\left| \sin(4\pi s) \right|$. In the third example, we consider a discontinuous profile with half of the data points uniformly distributed on the plane $x=1$ given by $h_+(y)=2 \sin(2\pi y)$ and the other half set of sampling points uniformly on $x=-1$ with $h_-(y)=-2\sin(2\pi y)$ with $y\in[-1,1]$.

Unless specified otherwise, we solve the advection equation (\ref{Eqn:advection}) with a reversible cosine velocity $c(s)=\cos\left( 2\pi s/T \right)$ and obtain the solutions at the final time $T=4$. For the semi-Lagrangian method, we solve the characteristics using the fourth-order Runge-Kutta method with a time step $10^{-3}$ to construct a flow map of size $\Delta t=10^{-1}$. Therefore, to obtain the solution at $T=4$, we iterate and interpolate the solution 40 times.



\begin{figure}[!h]
\centering
\includegraphics[width=0.95\textwidth]{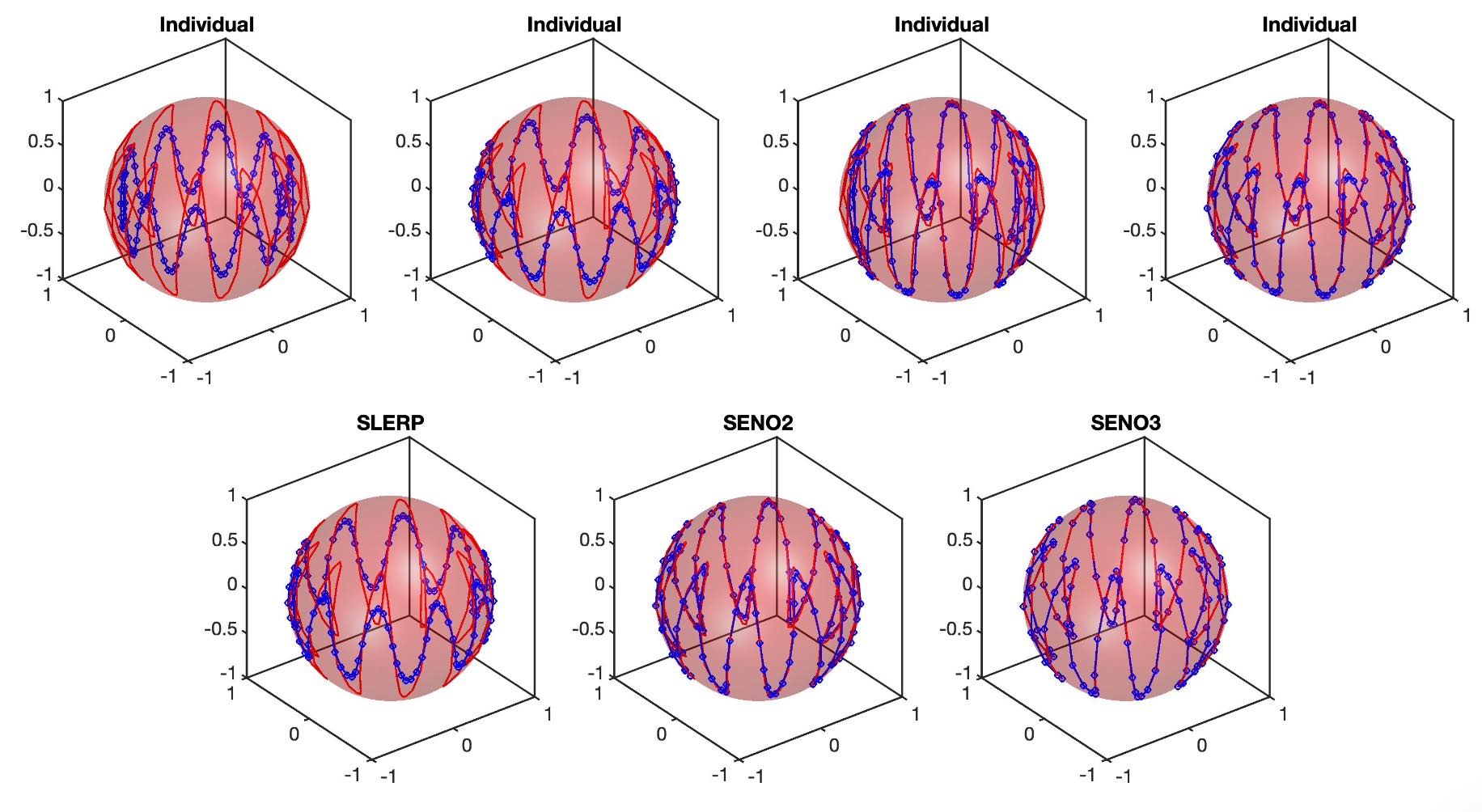}
\caption{(Section \ref{SubSec:ExSmooth} with the reversible cosine velocity) The exact solution at the final time is plotted using red solid line. (First row) The computed solutions using the component-by-component approach with linear interpolation without and with projection, \textsf{pchip} without and with projection. (Second row) The computed solutions using SLERP, SENO2 and SENO3.} 
\label{Fig:202}
\end{figure}


\begin{figure}[!h]
\centering
(a)\includegraphics[trim=0 0 10 0, clip, width=0.45\textwidth]{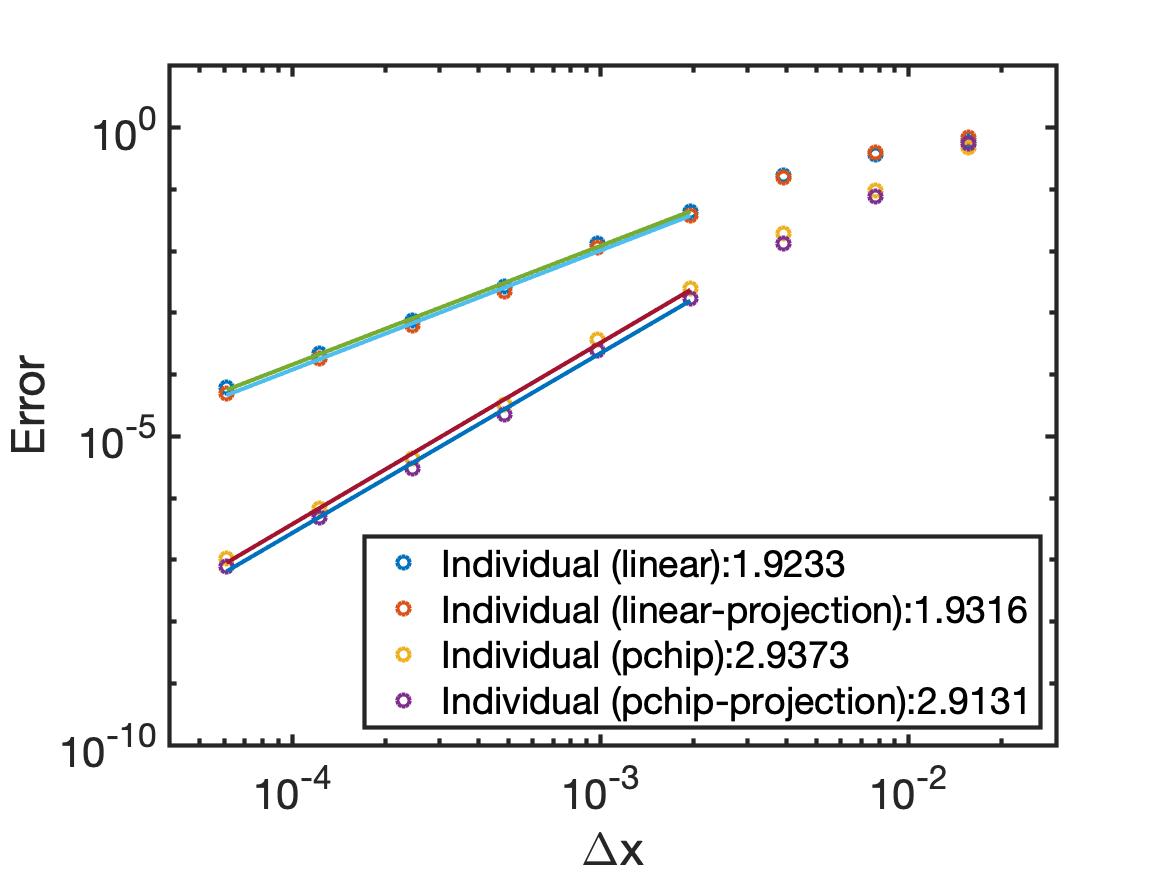} 
(b)\includegraphics[trim=0 0 10 0, clip, width=0.45\textwidth]{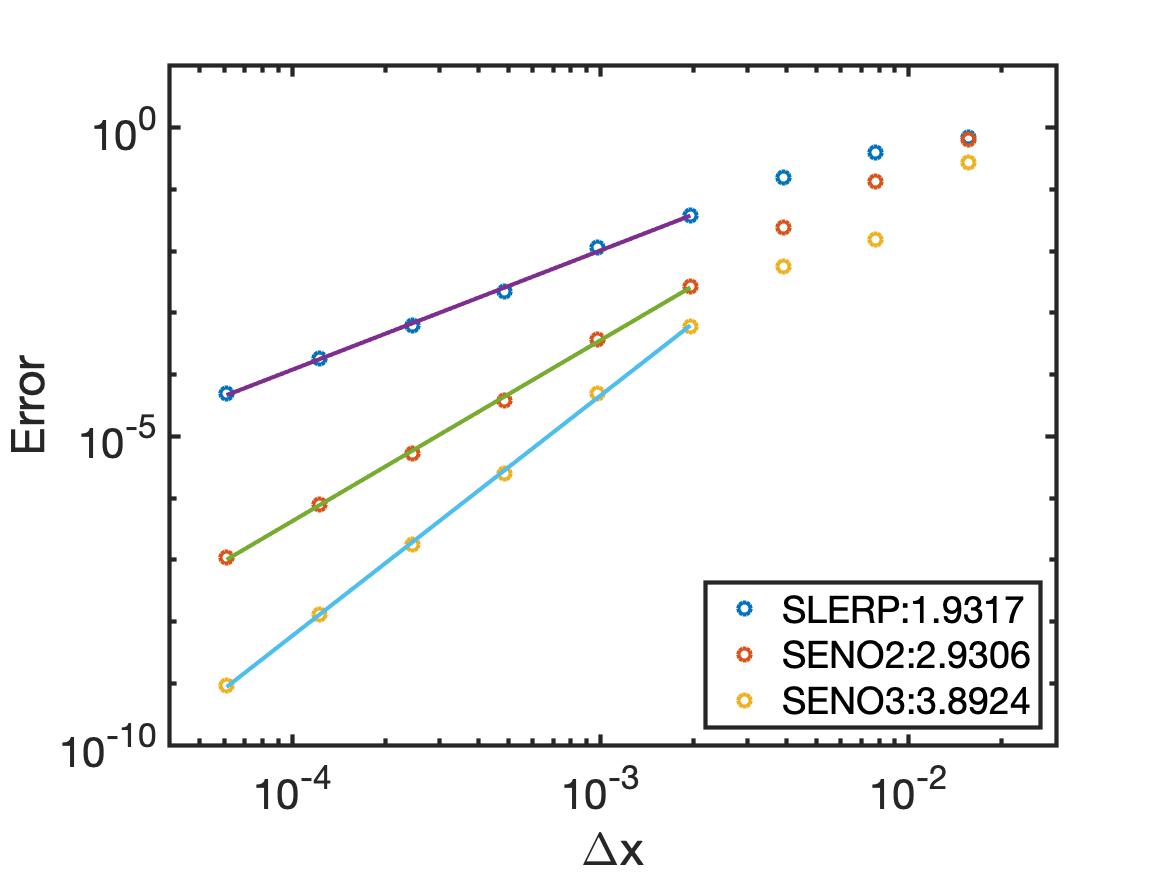} \\
(c)\includegraphics[trim=0 0 10 0, clip, width=0.45\textwidth]{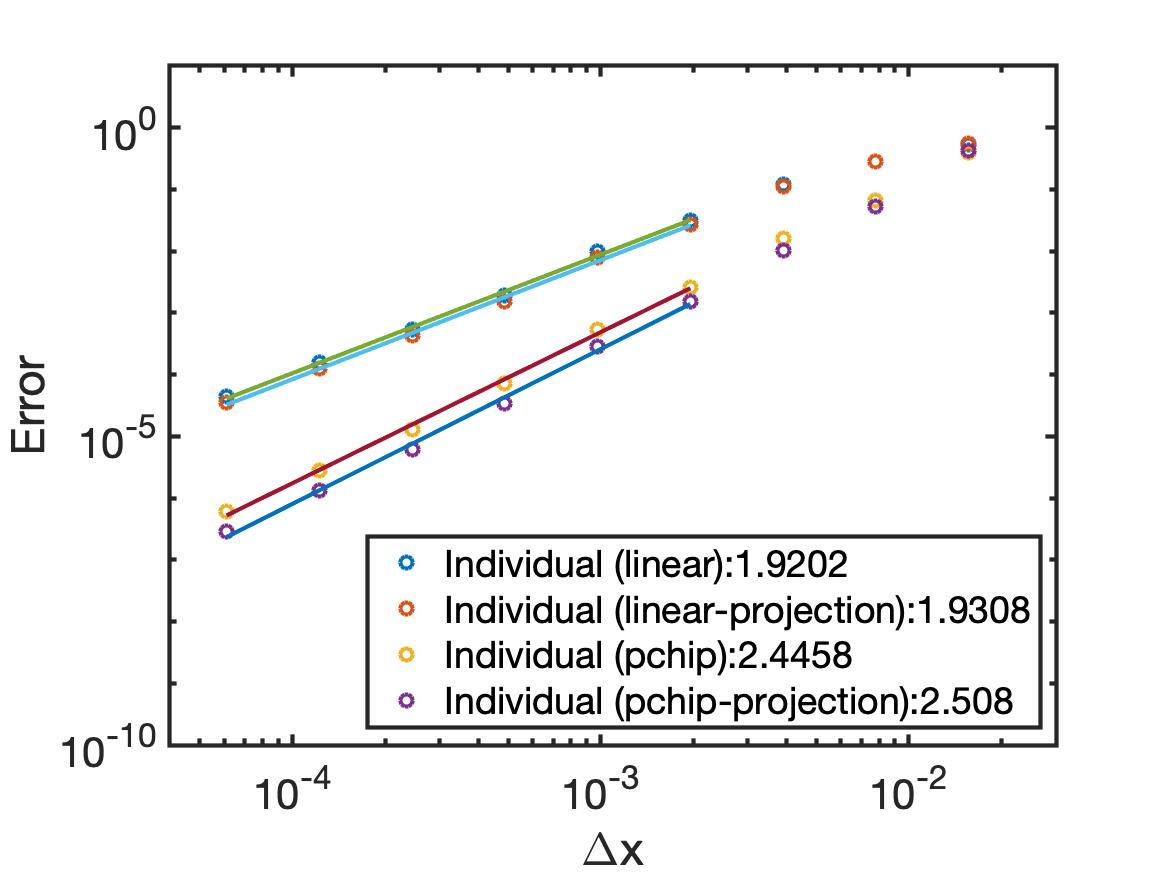}
(d)\includegraphics[trim=0 0 10 0, clip, width=0.45\textwidth]{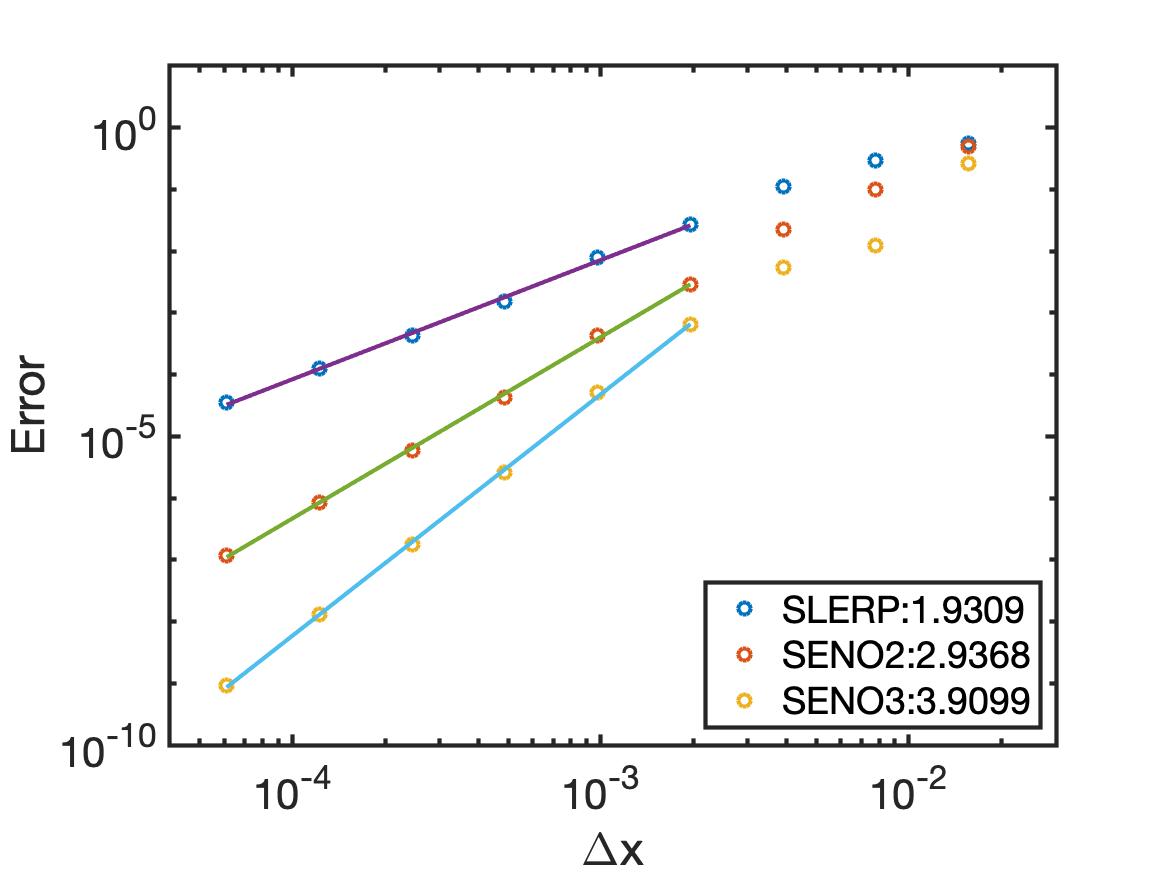}
\caption{(Section \ref{SubSec:ExSmooth} with the reversible cosine velocity) The $L_1$-error of the solutions obtained by (a) solving the PDE corresponding to individual components and (b) applying SLERP, SENO2 and SENO3. The $L_2$-error of the solutions obtained by (c) solving the PDE corresponding to individual components and (d) applying SLERP, SENO2 and SENO3.} 
\label{Fig:202_Error}
\end{figure}

\subsection{A Smooth Initial Condition}
\label{SubSec:ExSmooth}

We consider the smooth initial condition (\ref{Eqn:ExProjection1}) given by $g(s)=g_1(s)$ in this first example. Figure \ref{Fig:202} show the numerical solution using various semi-Lagrangian approaches computed on a relatively coarse mesh with $N=128$ on the interval $s=[0,1]$ with the velocity given by the reversible cosine velocity $c(s)$, respectively. We, in the first row of these figures, consider some component-by-component methods given by linear interpolation (which interpolates the $x$-, $y$-, and $z$-components individually), linear interpolation with projection (i.e., with a projection back onto the unit sphere after each interpolation), the function \textsf{pchip} implemented in \textsf{MATLAB}, and also \textsf{pchip} with projecting the solution back to $\mathbb{S}^2$. In the second row of these figures, we plot the solutions computed using the SLERP and our proposed SENO2 and SENO3 interpolations.

Because both linear interpolation and SLERP are developed based on some linear reconstructions, the solutions are of lower order. We see that the schemes are highly dissipative and cannot maintain the extrema in the solution profile. For the linear interpolation method without projection back to the sphere, we also see that the solution curve does not stay on the unit sphere $\p \notin \mathbb{S}^2$ but shrinks toward the origin. We can observe similar behavior for the solution obtained by the high-order shape-preserving piecewise cubic interpolation (i.e., \textsf{pchip}), for example as shown in the third subfigure on the first row in Figure \ref{Fig:202}. Since there is no control over the length of the solution $\|\p\|$, a projection is necessary to constrain that $\p \in\mathbb{S}^2$.

To confirm the accuracy of the proposed numerical schemes, we also perform a convergence test for both velocity fields. Since the exact solution to these test examples is the same as the initial condition, we define the following $L_1$- and $L_2$-errors given by
$E_1 = \int_0^1 \|\p(s)-\p_0(s)\|_1 \, ds$ and $E_2 = \left[\int_0^1 \|\p(s)-\p_0(s)\|_2^2 \, ds \right]^{1/2}$
where the initial condition is given by $\p_0(s_i)=(x_i^0,y_i^0,z_i^0)$. 

We collect the errors corresponding to the mesh numbers ranging from $N=64$ to 16384 and show them in Figure \ref{Fig:202_Error}. The $L_1$-errors are shown in (a-b), while the $L_2$-errors are plotted in (c-d). We see that the linear interpolation methods give second-order accurate solutions, while the \textsf{MATLAB} \textsf{pchip} interpolation provides a $L_1$ third-order accurate solution (but is dropped to somewhere between two and three when we measure the error in $L_2$-norm). And these orders are independent of the projection operator. Although there is no guarantee that the numerical solution stays on $\mathbb{S}^2$, the convergence matches those without an extra projection step. The accuracy in the solution obtained using SLERP is similar to that of linear interpolation since both methods use linear reconstruction. Therefore, we also see a second-order convergence in the SLERP method. As demonstrated in \cite{fonleu23}, both SENO2 and SENO3 are high-order interpolating schemes, and we observe a third- and fourth-order convergence in the numerical solutions.

\begin{figure}[!h]
\centering
\includegraphics[trim=60 0 100 0, clip, width=0.32\textwidth]{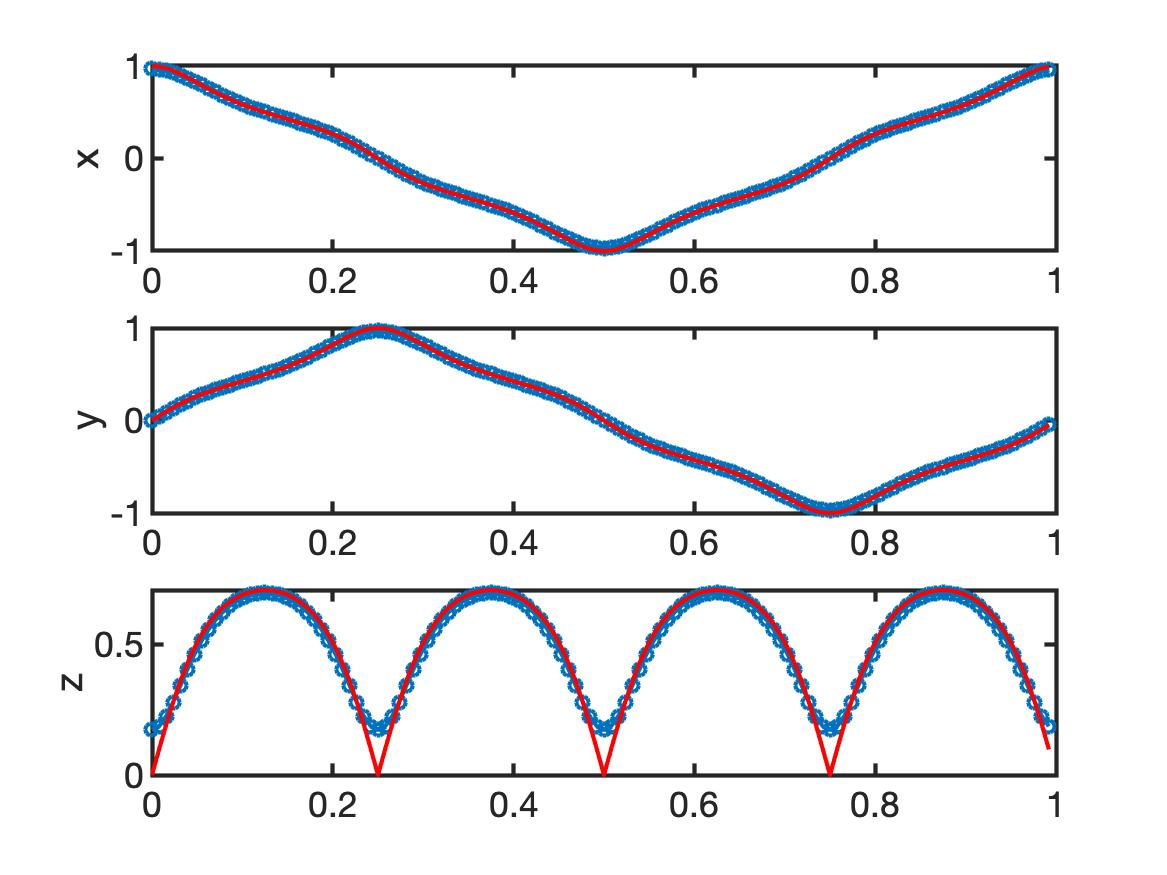} 
\includegraphics[trim=60 0 100 0, clip, width=0.32\textwidth]{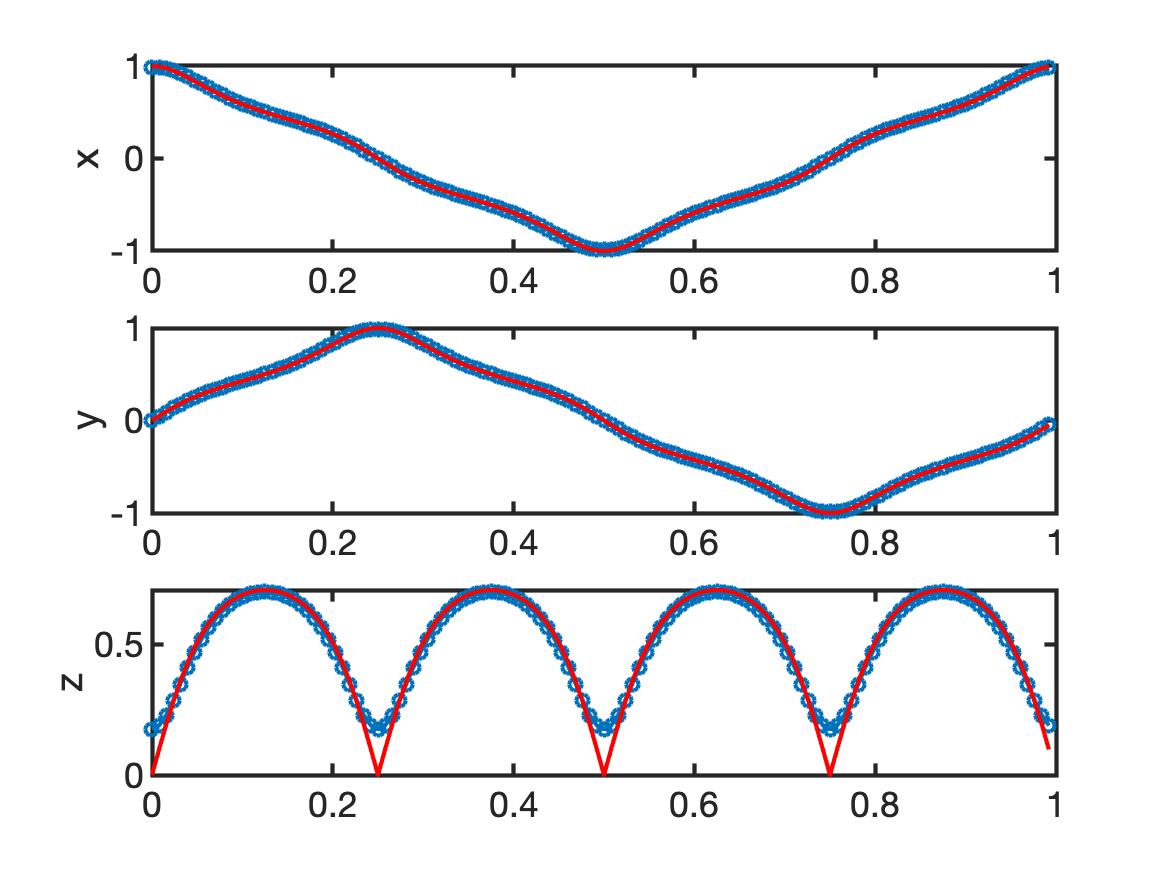} \\
\includegraphics[trim=60 0 100 0, clip, width=0.32\textwidth]{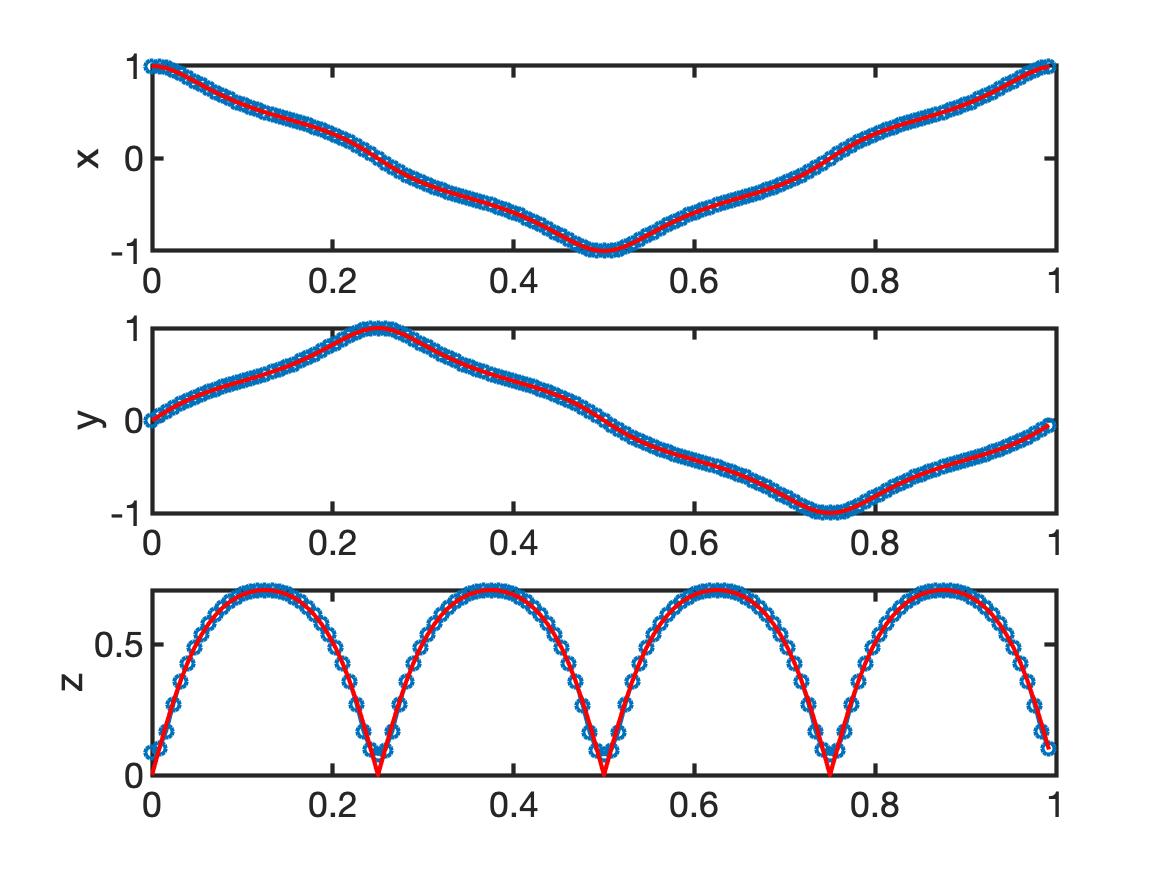} 
\includegraphics[trim=60 0 100 0, clip, width=0.32\textwidth]{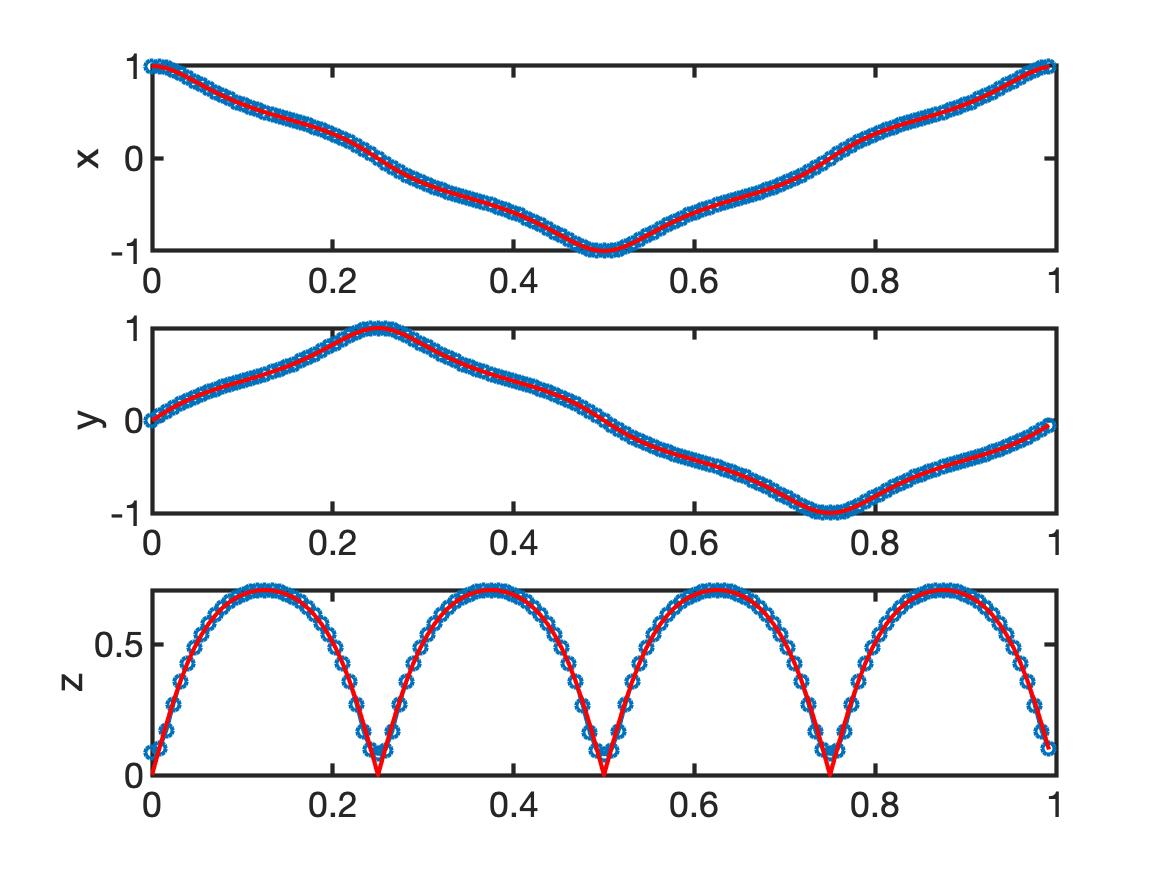} \\
\includegraphics[trim=60 0 100 0, clip, width=0.32\textwidth]{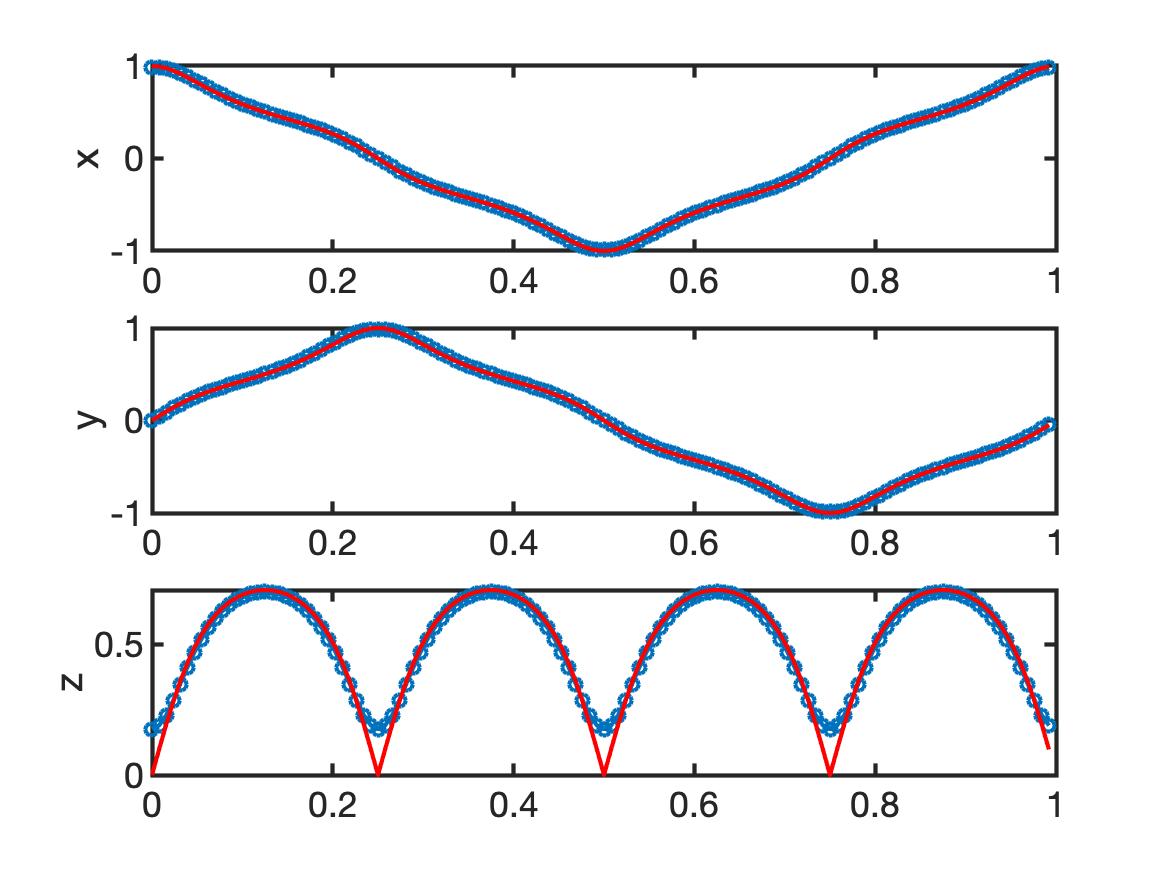} 
\includegraphics[trim=60 0 100 0, clip, width=0.32\textwidth]{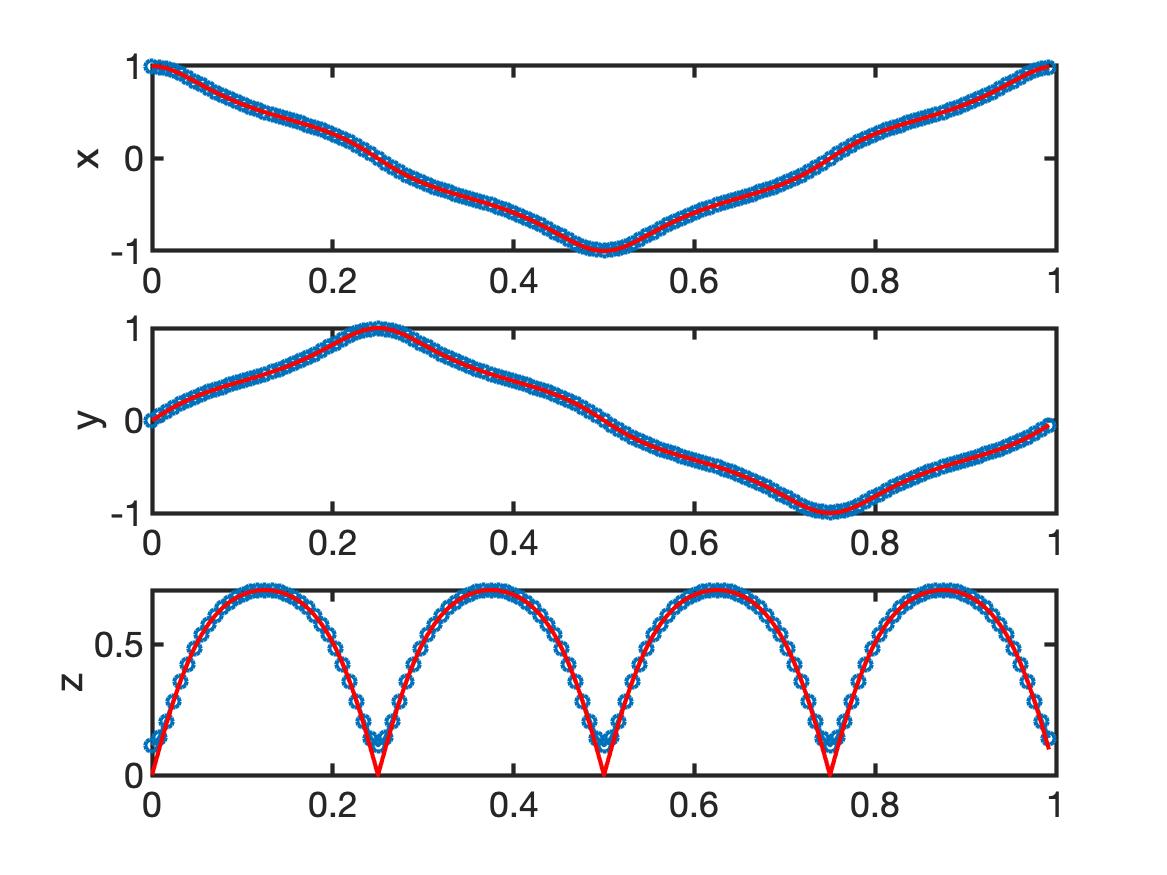} 
\includegraphics[trim=60 0 100 0, clip, width=0.32\textwidth]{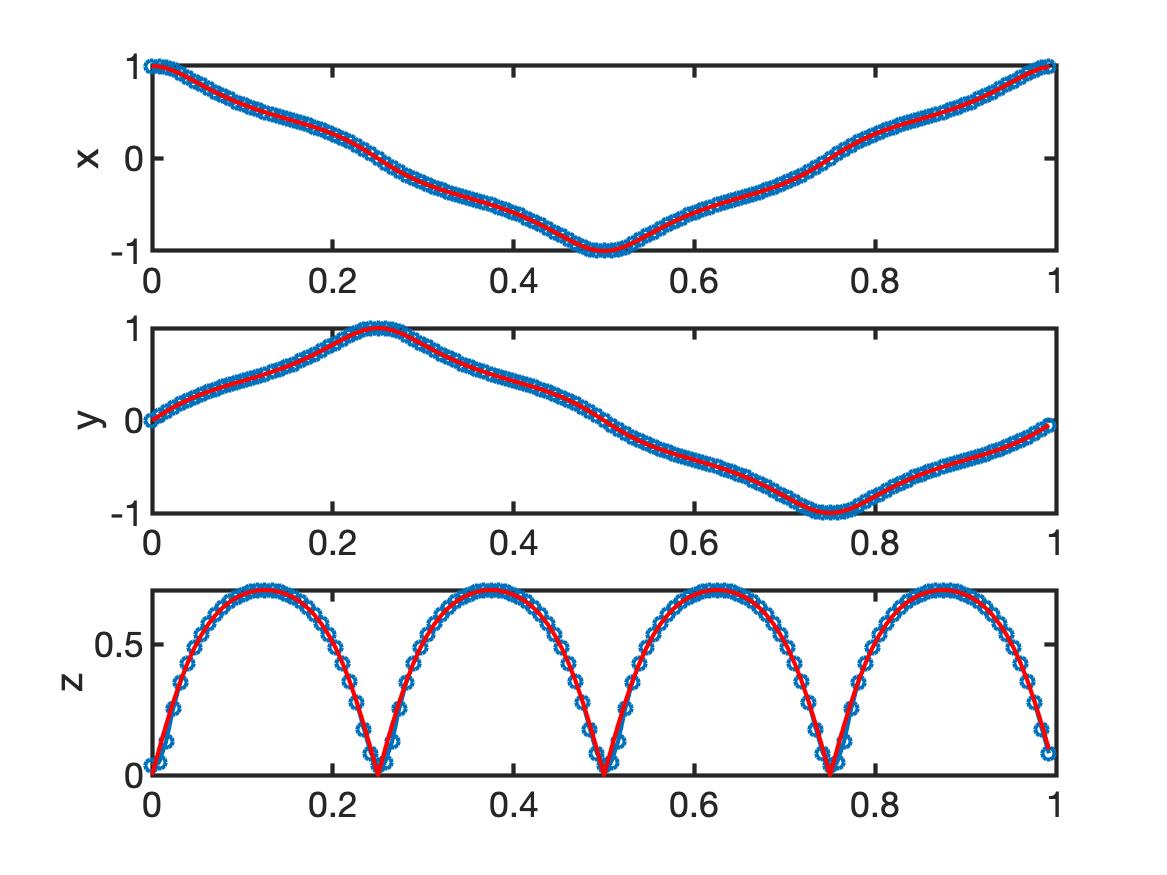} 
\caption{(Section \ref{SubSec:Kinks} with the reversible cosine velocity) The exact solution at the final time is plotted using red solid line.  (First row) The computed solutions using the component-by-component approach with linear interpolation without and with projection. (Second row) The computed solutions using the component-by-component approach with \textsf{pchip} without and with projection. (Third row) The computed solutions using SLERP, SENO2 and SENO3.} 
\label{Fig:302_128_components}
\end{figure}

\begin{figure}[!h]
\centering
\includegraphics[width=0.95\textwidth]{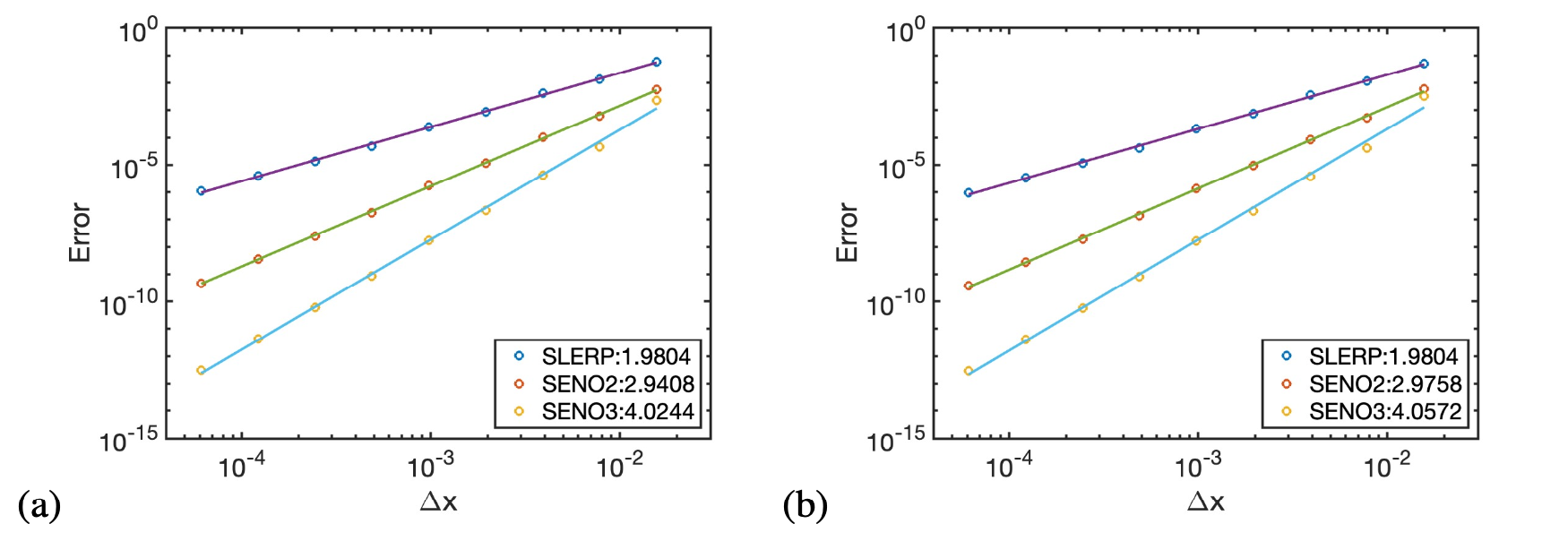}
\caption{{(Section \ref{SubSec:Kinks} with the reversible cosine velocity) The (a) $L_1$-error and (b) $L_2$-error of the solutions obtained by applying SLERP, SENO2 and SENO3.}} 
\label{Fig:302_Error}
\end{figure}

\subsection{An Initial Condition with Kinks}
\label{SubSec:Kinks}

In this example, we consider the initial condition with four kinks given by $g_2(s)$. Figure \ref{Fig:302_128_components} shows the solutions computed using various algorithms at the final time $T=4$ using the coarse mesh $N=128$. Since the size of the periodic domain is one, the motion has completed four periods. The first two rows in this figure show the solutions obtained by individually treating the advection equation component and see that the solutions do not differ much on whether we have the extra projection step. However, because of the extra numerical dissipations in the linear interpolation, the solutions are much smoother, and the method cannot preserve the kink in the $z$-component. While the \textsf{pchip} uses cubic reconstruction, the solution is more accurate and better resolves the kink. However, we observe some asymmetry in the solution. The solution tends to flatter the right-hand side after the kink in the $z$-component of $\p$. In the last row of Figure \ref{Fig:302_128_components}, we show the computed solutions using SLERP, SENO2, and SENO3. These solutions are all symmetric near the kinks, and we can much better resolve these singularities using SENO3. {Figure \ref{Fig:302_Error} illustrates the $L_1$- and $L_2$-errors in the solution, specifically in the region outside of the kink characterized by the $z$-component of the exact solution that exceeds 0.5. We also see that the numerical accuracy matches the expected values very well.}

\begin{figure}[!h]
\centering
\includegraphics[width=0.95\textwidth]{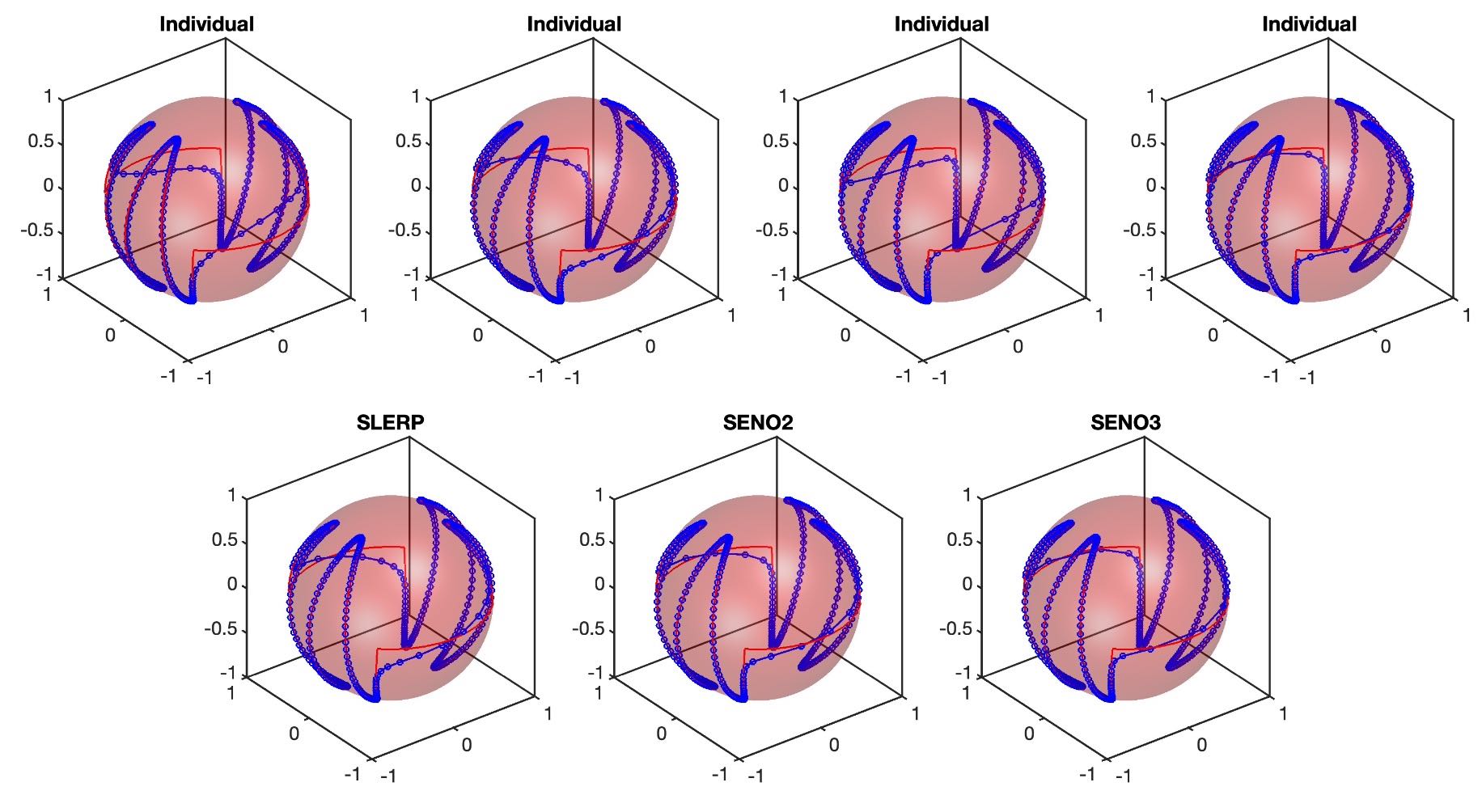}
\caption{(Section \ref{SubSec:ExDiscontinuous} with the reversible cosine velocity) The exact solution at the final time is plotted using red solid line. (First row) The computed solutions with $N=512$ using the component-by-component approach with linear interpolation without and with projection, \textsf{pchip} without and with projection. 
(Second row) The computed solutions using SLERP, SENO2 and SENO3.} 
\label{Fig:402_512}
\end{figure}


\begin{figure}[!h]
\centering
\includegraphics[trim=60 0 100 0, clip, width=0.32\textwidth]{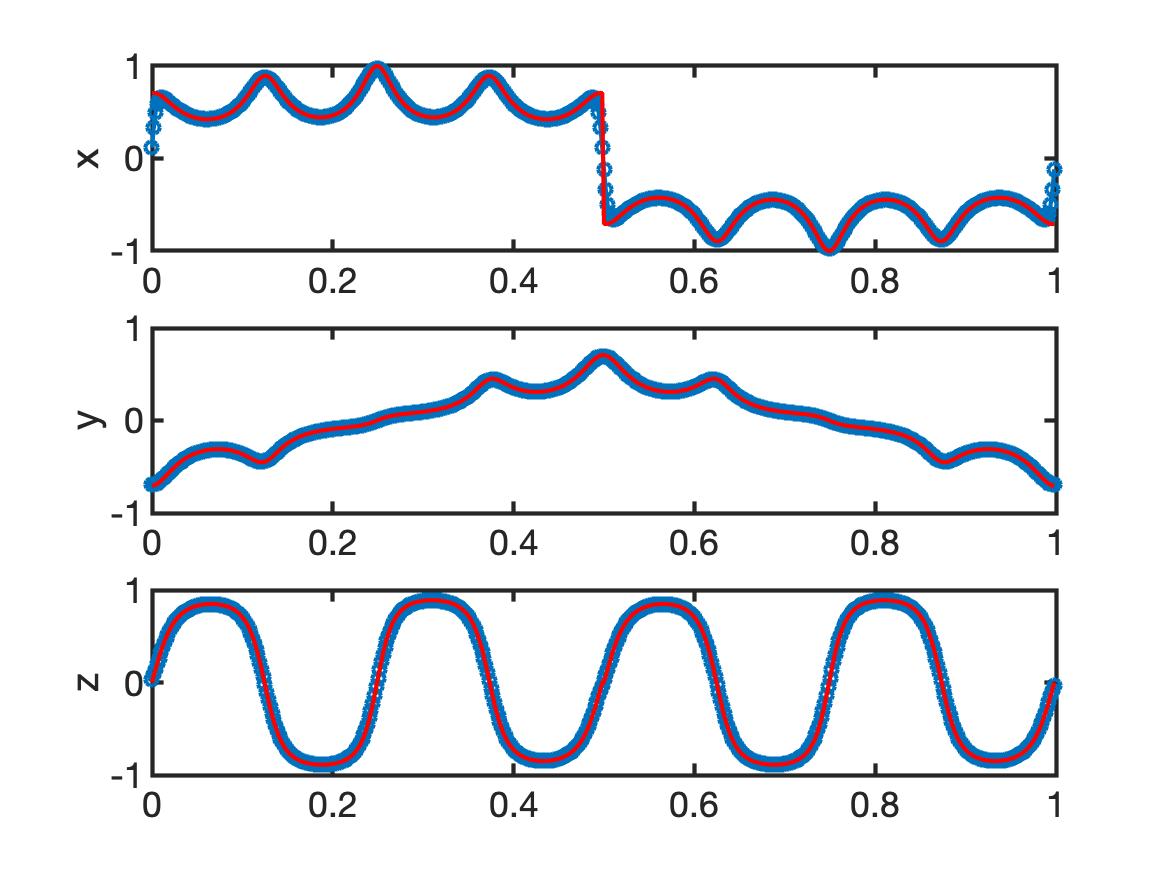} 
\includegraphics[trim=60 0 100 0, clip, width=0.32\textwidth]{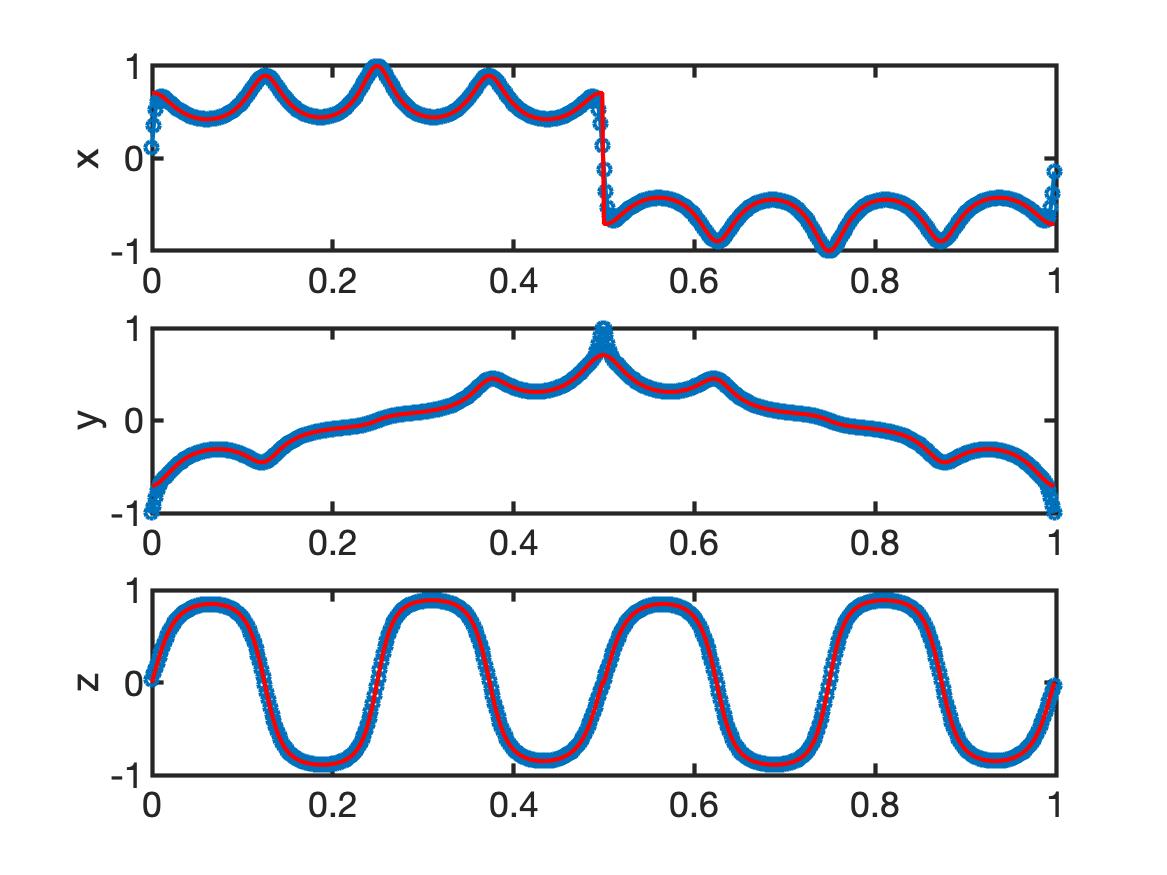} \\
\includegraphics[trim=60 0 100 0, clip, width=0.32\textwidth]{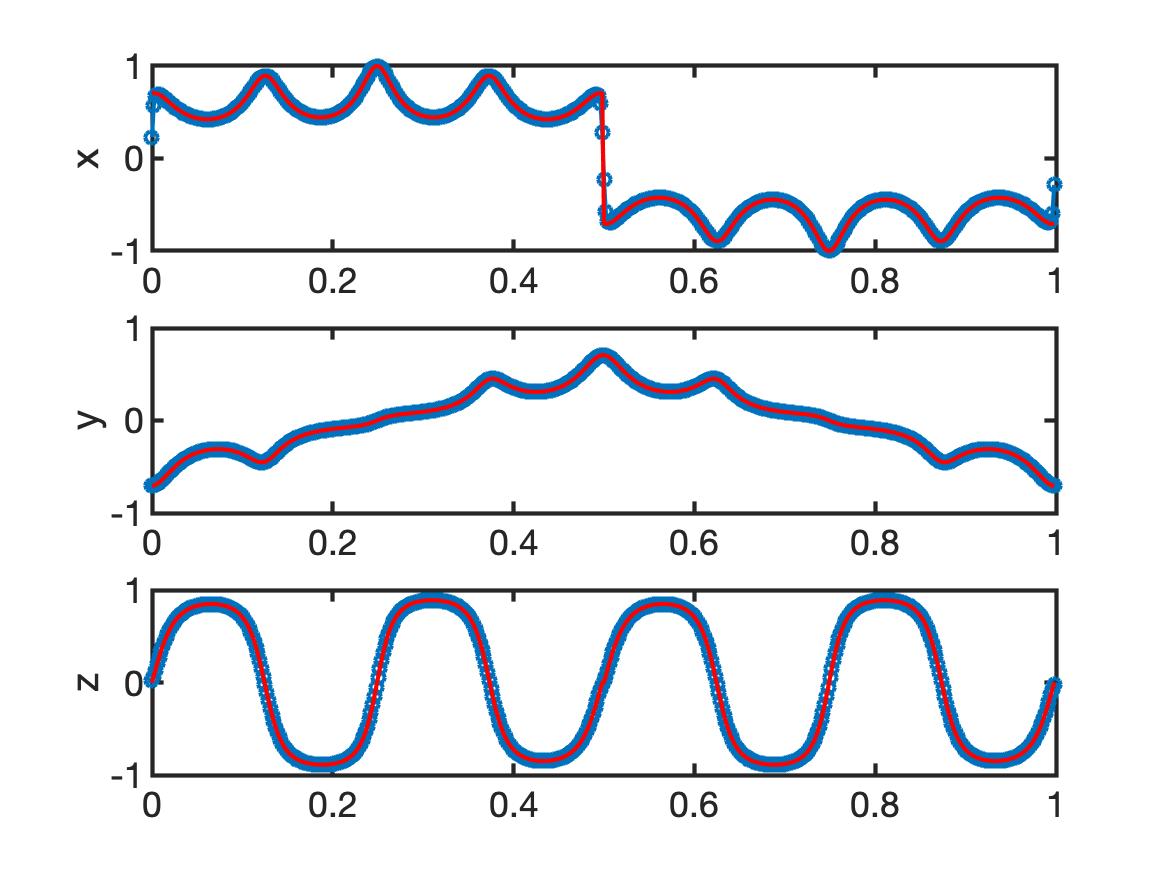} 
\includegraphics[trim=60 0 100 0, clip, width=0.32\textwidth]{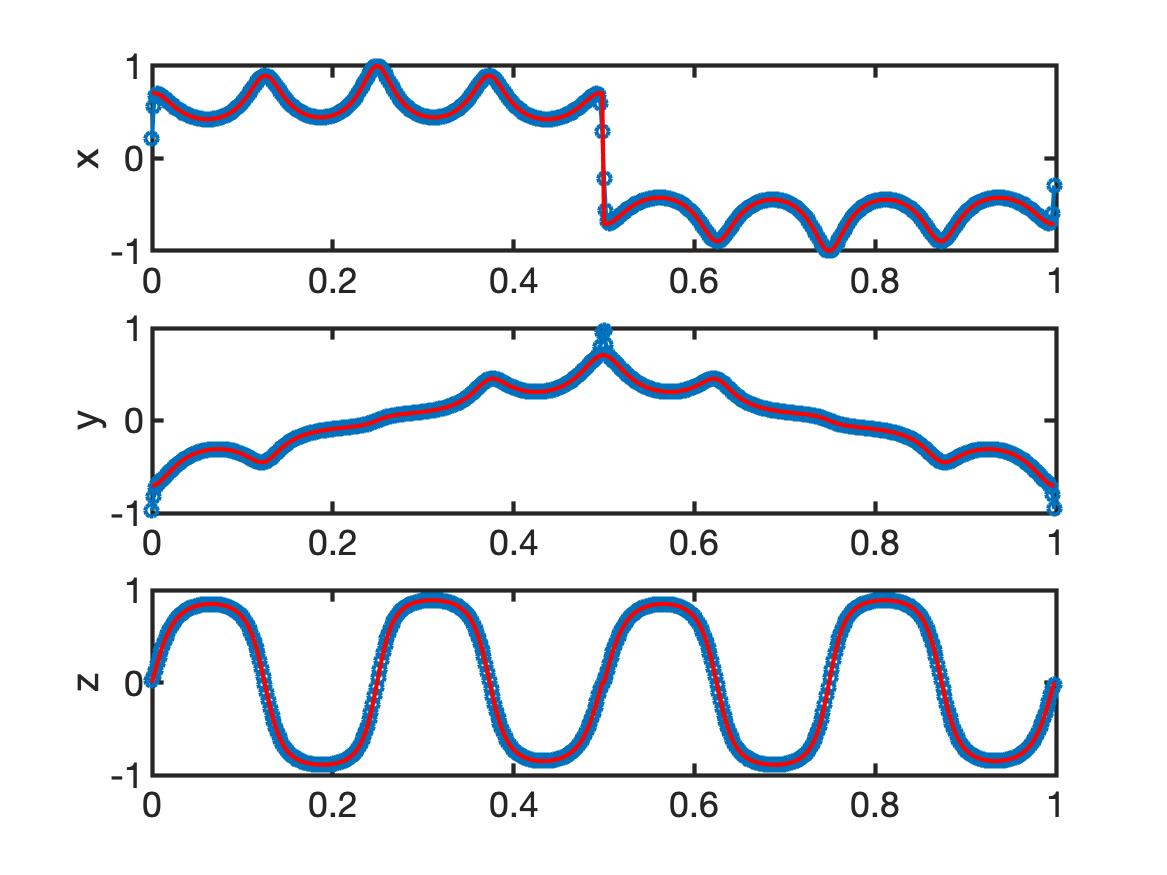} \\
\includegraphics[trim=60 0 100 0, clip, width=0.32\textwidth]{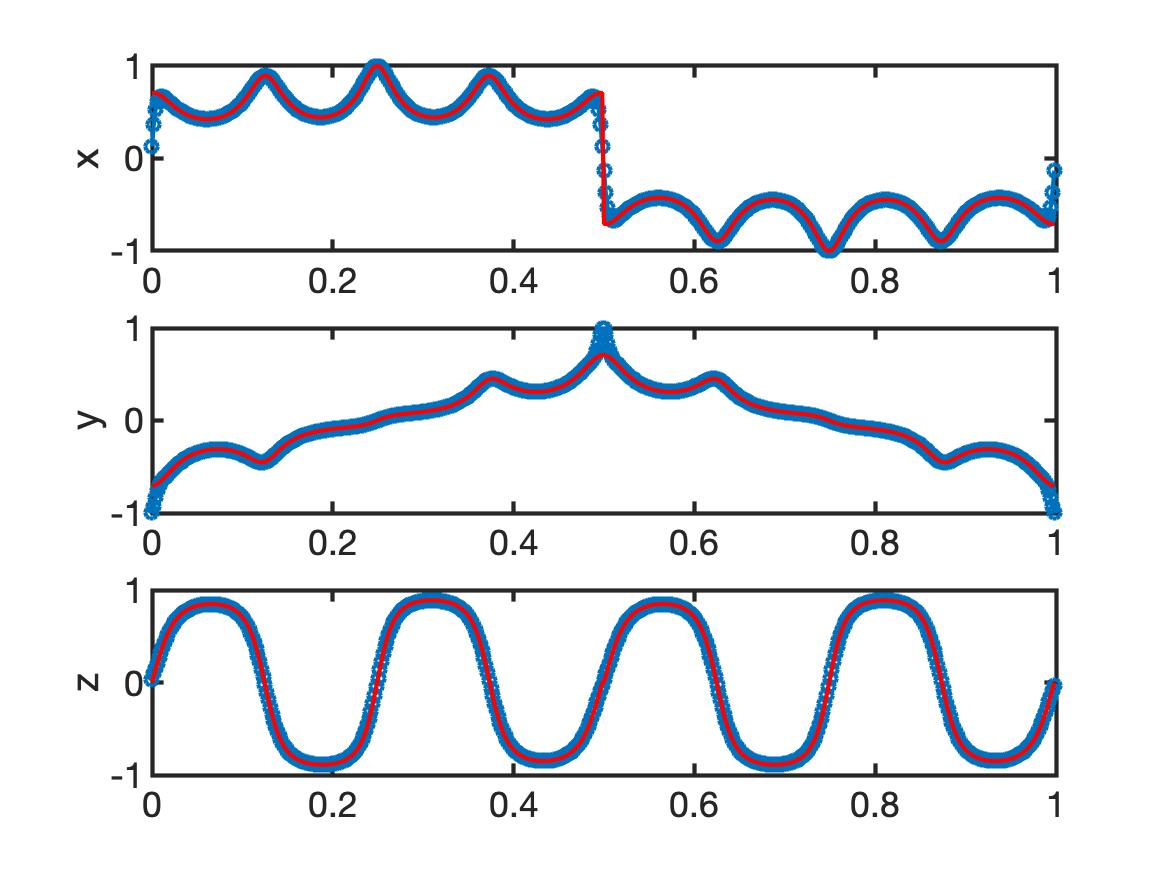} 
\includegraphics[trim=60 0 100 0, clip, width=0.32\textwidth]{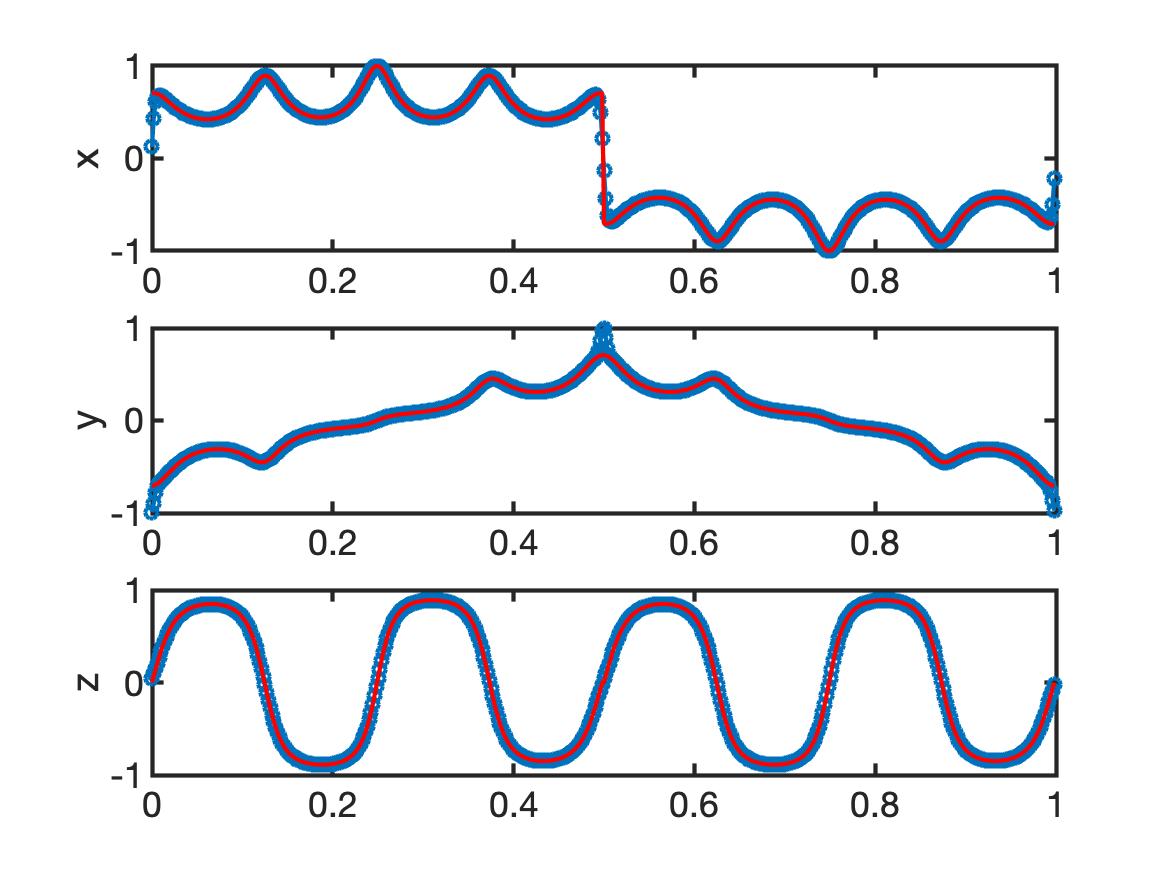} 
\includegraphics[trim=60 0 100 0, clip, width=0.32\textwidth]{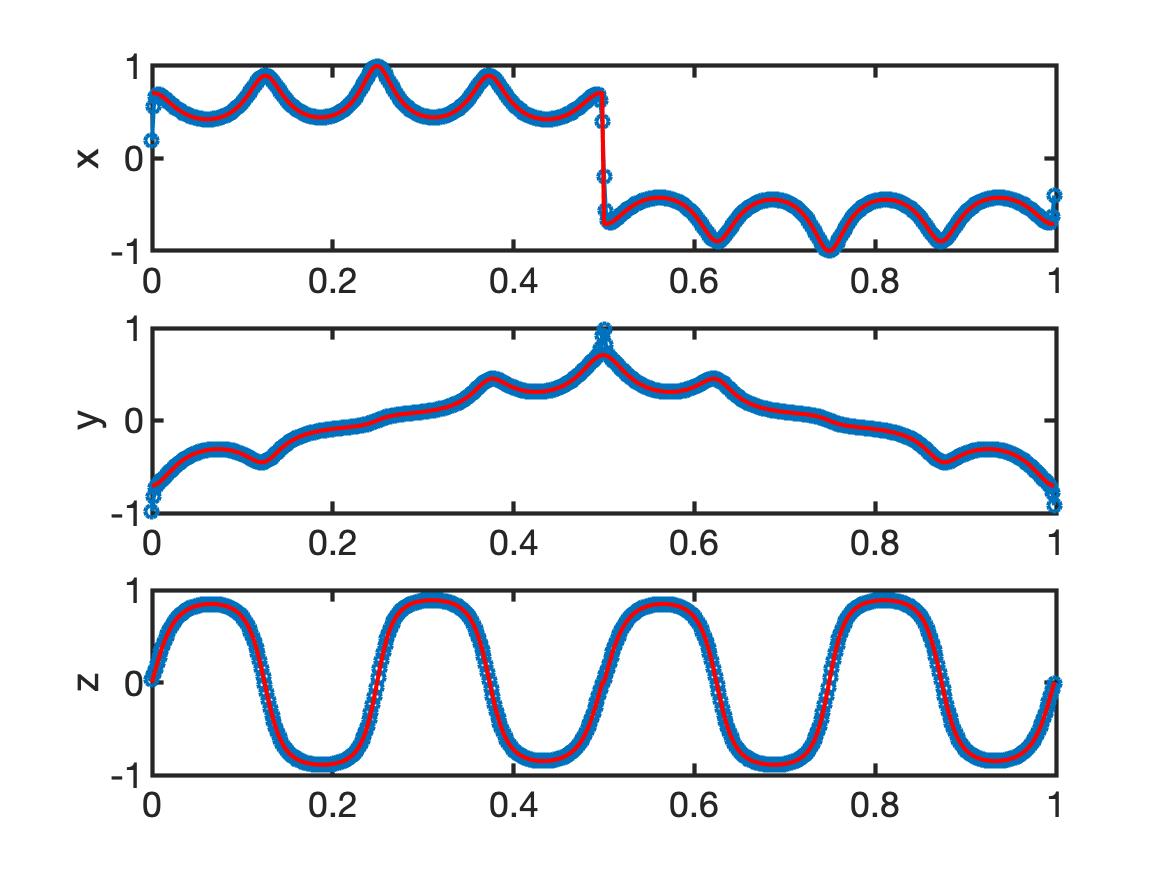} 
\caption{(Section \ref{SubSec:ExDiscontinuous} with the reversible cosine velocity) The exact solution at the final time is plotted using red solid line.  (First row) The computed solutions with $N=512$ using the component-by-component approach with linear interpolation without and with projection. (Second row) The computed solutions using the component-by-component approach with \textsf{pchip} without and with projection. (Third row) The computed solutions using SLERP, SENO2 and SENO3.} 
\label{Fig:402_512_components}
\end{figure}

\subsection{A Discontinuous Initial Condition}
\label{SubSec:ExDiscontinuous}

This example considers the evolution of an initial discontinuous profile by constructing the functions $h_{\pm}(y)$ and projecting them onto the unit sphere. Such an initialization gives two jumps in the $x$-component of the function $f$. Figures \ref{Fig:402_512} shows the solutions computed on a coarse mesh $N=512$ under the constant velocity motion. Once again, we are looking at the solutions at the final time $T=4$ so that these solutions have completed four periods. To better visualize the discontinuities in the solution, we plot each component of $\p$ in Figure \ref{Fig:402_512_components}. We see that SENO3 provides the most accurate solution to resolving the shock using only one to two grid points.

\section{Conclusion}

{In conclusion, we have successfully developed a numerical scheme for solving the advection equation of $\mathbb{S}^2$-valued functions of real variables. This scheme allows for high-order modeling of the time-evolution of a $\mathbb{S}^2$-valued mapping on the real line. Our approach builds upon the semi-Lagrangian method for the linear scalar advection equation and also the SENO interpolation method to address the challenges posed by $\mathbb{S}^2$-functions with kinks or sharp discontinuities in their components. Future work includes an extension to advection equations of SO(3), i.e., real $3\times 3$ orthogonal matrices of unit determinant.}

\section*{Acknowledgment}
The work of Leung was supported in part by the Hong Kong RGC grant 16302223.

\bibliographystyle{plain}
\bibliography{syleung}


\end{document}